\documentclass[reqno, 11pt]{amsart}
\title{Quenched invariance principle for random walks on Delaunay triangulations}
\author[A. Rousselle]{Arnaud Rousselle$^\ddagger$
}
\thanks{
$\ddagger$ Normandie Universit\'e, Universit\'e de Rouen, 
Laboratoire de Math\'ematiques Rapha\"el Salem,
CNRS, UMR 6085, Avenue de l'universit\'e, BP 12, 
76801 Saint-Etienne du Rouvray Cedex, 
France.
\newline
and
\newline
Laboratoire Modal'X, Universit\'e Paris Ouest Nanterre La D\'efense, 200 avenue de la R\'epublique, 92000 Nanterre,
France
\newline
E-mail address: arnaud.rousselle@u-paris10.fr}

\usepackage{amssymb}
\usepackage{stmaryrd}
\usepackage{latexsym}
\usepackage{mathtools}
\usepackage{epsfig}
\usepackage{epsf}
\usepackage{float}
\usepackage{a4wide}
\usepackage{hyperref}
\usepackage{latexsym,amssymb,amsmath,epsfig,psfrag}
\usepackage{alltt}
\usepackage{cite}
\usepackage{color}
\usepackage{wasysym}
\def\0{\mathbf{0}}
\def\x{\mathbf{x}}
\def\d{\mathrm{d}}

\def\z{\mathbf{z}}
\def\N{\mathcal{N}}

\def\P{\mathcal{P}}
\def\B{\mathcal{B}}

\def\A{\mathbf{A}}

\def\E{\mathcal{E}}

\def\x{\mathbf{x}}
\def\y{\mathbf{y}}
\def\z{\mathbf{z}}

\def\RR{\mathbb{R}}

\def\ZZ{\mathbb{Z}}

\def\NN{\mathbb{N}}

\def\xiz{{\xi^0}}

\def\A{\mathcal{A}}

\def\Vor{\operatorname{Vor}}

\newcommand{\sfrac}[2]{\kern.1em
        \raise.5ex\hbox{$#1$}\kern-.1em
        /\kern-.15em\lower.25ex\hbox{$#2$}}
\newtheorem*{ass*}{Assumptions}
\newtheorem{defi}{Definition}
\newtheorem{lemm}[defi]{Lemma}

\newtheorem{prop}[defi]{Proposition}
\newtheorem{coro}[defi]{Corollary}
\newtheorem{theo}[defi]{Theorem}
\newtheorem{exem}[defi]{Example}
\newtheorem{rem}[defi]{Remark}
\newtheorem{fact}[defi]{Fact}

\newenvironment{dem}{\noindent {\bf Proof.}}
                    {\hfill $\square$}
\newenvironment{demof}[1]{\noindent{\bf Proof of #1.}}
                    {\hfill $\square$}

\renewcommand{\div}{\operatorname{div}}
\def\bex{\begin{exem} \em }
\def\eex{\end{exem} }
\def\brem{\begin{rem} \em}
\def\erem{\end{rem} }

\def\deg{\operatorname{deg}}
\def\DT{\operatorname{DT}}

\def\diam{\operatorname{diam}}
\begin{document}
\begin{abstract}
We consider simple random walks on Delaunay triangulations generated by point processes in $\mathbb{R}^d$. Under suitable assumptions on the point processes, we show that the random walk satisfies an almost sure (or quenched) invariance principle. This invariance principle holds for point processes which have clustering or repulsiveness properties including Poisson point processes,  Mat\'ern cluster and Mat\'ern hardcore processes. The method relies on the decomposition of the process into a martingale part and a \emph{corrector} which is proved to be negligible at the diffusive scale. 
\end{abstract}
\date{\today}
\maketitle



\maketitle

\section{Introduction}
Let us first describe the model. Given an infinite locally finite subset $\xi$  of $\RR^d$, the Voronoi tessellation associated with $\xi$ is the collection of the Voronoi cells: $$\mathrm{Vor}_\xi(\x):=\{x\in\RR^d\, :\,\Vert x-\x\Vert\leq\Vert x-\y\Vert,\forall\y\in\xi \},\quad\x\in\xi.$$
The point $\x$ is called the nucleus or the seed of the cell. The Delaunay triangulation $\operatorname{DT}(\xi)$ of $\xi$ is the dual graph of its Voronoi tiling. It has $\xi$ as vertex set and there is an edge between $\x$ and $\y$ in $\operatorname{DT}(\xi)$ if $\mathrm{Vor}_\xi(\x)$ and $\mathrm{Vor}_\xi(\y)$ share a $(d-1)$-dimensional face. Another useful characterization of $\operatorname{DT}(\xi)$ is the following: a simplex $\Delta$ is a cell of $\operatorname{DT}(\xi)$ \emph{iff} its circumscribed sphere has no point of $\xi$ in its interior. Recall that this is a well defined triangulation when $\xi$ is in general position. In the sequel, we denote by $\N$ (resp. $\N_0$) the set of infinite locally finite subsets of $\RR^d$ (resp. the set of  infinite locally finite subsets of $\RR^d$ containing 0).

Given a realization $\xi$ of a suitable point process with law $\P$, we consider the variable speed nearest-neighbor random walk $(X_t)_{t\geq 0}$ on the Delaunay triangulation of $\xi$, that is the Markov process with generator:
\begin{equation}\label{EqDefGene}
\mathcal{L}^\xi f(x):=\sum_{y\in\xi}c^\xi_{x,y}\left(f(y)-f(x)\right),\qquad x\in\xi,
\end{equation} 
where $c^\xi_{x,y}$ is the indicator function of `$y\sim x$ in $\DT(\xi)$'. We denote by $P^\xi_x$ the law of the walk starting from $x\in\xi$ and by $E^\xi_x$ the corresponding expectation. We study for almost every realization of the point process the behavior of the random walk at the diffusive scale and we prove the following theorem.
\begin{theo}\label{QIP}
Assume that $\xi$ is distributed according to a simple, stationary, isotropic point process with law $\P$ a.s. in general position, satisfying assumptions {\bf (V)}, {\bf (SD)}, {\bf (Er)} and {\bf (PM)} (see Subsection \ref{SectCondPP}).

 For $\mathcal{P}-$a.e. $\xi$, for all $x\in\xi$, under $P^\xi_x$, the rescaled process $(X^\varepsilon_t)_{t\geq 0}=(\varepsilon X_{\varepsilon^{-2}t})_{t\geq 0}$ converges in law as $\varepsilon$ tends to $0$ to a Brownian motion with covariance matrix $\sigma^2I$ where $\sigma^2=\sigma^2_{\text{VSRW}}$ is positive and does not depend on $\xi$. 
\end{theo} 
Note that the assumptions of Theorem \ref{QIP} are satisfied by Poisson point processes, Mat\'ern hardcore processes and Mat\'ern cluster processes. The same result holds for the discrete-time nearest-neighbor random walk $(X_n)_{n\in \NN}$ with diffusion coefficient $\sigma^2_{\text{DTRW}}$ related with $\sigma^2_{\text{VSRW}}$ by the formula: $$\sigma^2_{\text{VSRW}}=\E_0\left[\deg_{\DT(\xi^0)}(0)\right]\sigma^2_{\text{DTRW}},$$ where $\E_0$ denotes the expectation w.r.t. the Palm measure $\P_0$ associated with the (stationary) point process with law $\P$. 

The main idea for proving such results is to show that the random walk behaves like a martingale up to a correction which is negligible at the diffusive scale. Actually, well-known arguments (see {\it e.g.} \cite[\S 3.3.1]{CFP}, \cite[p. 1340-1341]{BP} or \cite[\S 6.1 and \S 6.2]{BergerBiskup}) show that the last claim follows from Theorem \ref{TheoSousLine}.
\begin{theo}\label{TheoSousLine}
Under the assumptions of Theorem \ref{QIP}, there exists a so-called \emph{corrector} $\chi : \N_0\times \RR^d\longmapsto \RR^d$ such that for:
\begin{enumerate}
\item \label{ThSousLin1}$\varphi(\xi^0,x):=x-\chi (\xi^0,x)$ is \emph{harmonic at 0} for $\P_0-a.e.~\xi^0$, {\it i.e.}:
\[\sum_{x\in\xi^0}c^{\xi^0}_{0,x}\Vert \varphi(\xi^0,x)\Vert<\infty \mbox{ and }\mathcal{L}^{\xi^0}\varphi(\xi^0,0)=0\mbox{ for } \P_0-a.e.~\xi^0;\]
\item $\chi$ is a.s. sublinear:
\[\frac{\underset{x\in\xi^0\cap [-n,n]^d}{\max}\left\Vert\chi \left(\xi^0,x\right)\right\Vert}{n}\xrightarrow[n\to\infty]{\P_0-a.s.} 0.\]
\end{enumerate}
\end{theo}
The arguments to deduce Theorem \ref{QIP} from Theorem \ref{TheoSousLine} are rather standard. We have chosen not to do develop the arguments which can be found in the references cited above and we only indicate the main lines of the proof in Section \ref{ProofTHPrincQIP}. 
Various methods to prove quenched invariance principle for random walks among random conductances on $\ZZ^d$ or related models were developed during the last ten years (see \cite{SidoravS,BergerBiskup,BP,MP,BZ08,MathieuQuenched,BiskupLN,GZ13,ZZ13,BD14}). Theorem \ref{TheoSousLine} is proved by adapting the approach developed in \cite{BP}. Actually, we first prove the sublinearity of the corrector restricted to a suitable subgraph of the Delaunay triangulation and extend it by harmonicity. Let us note that this method was also successfully used in the context of random walks on complete graphs generated by point processes with jump probability which is a decreasing function of the distance between points in \cite{CFP}.

Recurrence and transience results for random walks on Delaunay triangulations generated by point processes were obtained in \cite{RTCRWRG} and an annealed invariance principle was proved in \cite{AIP}.
In \cite{FGG}, the existence of an harmonic corrector was recently established by a different (constructive) method. Nevertheless, the authors of this paper obtained the sublinearity of the corrector only in dimension 2. In order to extend the quenched invariance principle in higher dimensions, they suggested to prove full heat-kernel bounds similar to the one obtained by Barlow in \cite{Barlow} in the setting of random walks on supercritical percolations clusters. This approach would require much more sophisticated arguments and a better control of the regularity of the full graph than the one used in the present paper. It is worth noting that, as in \cite{RTCRWRG, AIP}, the arguments given in the sequel can be used to obtain quenched invariance principles for random walks on other graphs constructed according to the geometry of a realization of a point process. 
\subsection{Conditions on the point process}\label{SectCondPP}
In this section, we list the assumptions on the point process needed to obtain the quenched invariance principle. 

We assume that $\xi$ is distributed according to a simple and stationary point process with law $\P$. In the sequel, we will denote by $\E$ the expectation with respect to $\P$. We suppose that $\P$ is isotropic and almost surely in general position (see \cite{Zessin}): there are no $d+1$ points (resp. $d+2$ points) in a $(d-1)$-dimensional affine subspace (resp. in a sphere). We also assume that the point process satisfies:
\begin{itemize}
\item[{\bf (V)}]
there exists a positive constant $c_1$ such that for $L$ large enough:
\begin{equation*}
\P\big[\#\big([0,L]^d\cap\xi\big)=0\big]\leq e^{-c_1L^d}.
\end{equation*}
\end{itemize}

In order to prove the sublinearity of the corrector, we need to restrict the study to a subgraph of the Delaunay triangulation of $\xi$ which has good regularity properties. To this end, we will define a notion of `good boxes' that in particular allows us to bound the maximal degree of vertices in an infinite subgraph of the Delaunay triangulation of $\xi$. Precise definitions and assumptions are given below. 
\subsubsection{Good boxes, good points and the stochastic domination assumption}\label{SectGoodBoxes4}
For $s\in\NN^*$, let us divide $\RR^d$ into boxes of side $K:=\lceil 3 \sqrt{d}\rceil s$:
\begin{equation*}\label{EqDefBox}
B_\z=B^K_\z:=K\z+\Big[-\frac{K}{2},\frac{K}{2}\Big]^d,\quad\z\in\ZZ^d.
\end{equation*}
Each box $B_\z$ is then subdivided into smaller sub-boxes $b^\z_i$, $i=1,\dots, \lceil 3 \sqrt{d}\rceil^d $ of side $s$. 

A box $B_\z$, $\z\in\ZZ^d$, is called ($\alpha$-)\emph{nice} if each sub-box $b^\z_i$ of side $s$ in $B_\z$ satisfies $1\leq \#(\xi\cap b^\z_i)\leq \alpha s^d$. A box $B_\z$ is then said to be ($\alpha$-)\emph{good} if $B_{\z'}$ is $\alpha$-nice for each $\z'\in\ZZ^d$ with $\Vert \z'-\z\Vert_\infty\leq 1$. Writing $p^\text{site}_c(\ZZ^d)$ for the critical probability for site percolation in $\ZZ^d$, the stochastic domination hypothesis is the following.

\begin{itemize}
\item[{\bf (SD)}] For any $p^\text{site}_c(\ZZ^d)<p<1$, if $\alpha$ and $s_0$ are well chosen, for any $s\geq s_0$ the process of good boxes $\mathbb{X}:=\{X_\z=\mathbf{1}_{B_\z\mbox{ is good}},\,\z\in\ZZ^d\}$ dominates an independent site percolation process $\mathbb{Y}:=\{Y_\z,\,\z\in\ZZ^d\}$ on $\ZZ^d$ with parameter $p$.
\end{itemize}

If {\bf (SD)} is satisfied, we can find a coupling $\mathcal{P}^{K,p}$ of the processes $\mathbb{X}$ and $\mathbb{Y}$ such that $B_\z$ is good when $Y_\z=1$. With a slight abuse of notation, we will omit the superscript and denote by $\mathcal{P}$ the probability measure $\mathcal{P}^{K,p}$ on $\widehat{\mathcal{N}}:=\mathcal{N}\times \{0,1\}^{\mathbb{Z}^d}$ whose marginal distributions are respectively the law of the point process $\xi$ and the law of the  independent site percolation process $\mathbb{Y}$ with parameter $p$ being fixed large enough. The precise value of $p$ is not stated explicitely but we assume that it is large enough to ensure that all the percolation results we need are satisfied. A generic element of $\widehat{\mathcal{N}}$ is denoted by $\widehat{\xi}=(\xi,(y_\z)_{\z\in\mathbb{Z}^d})$. 

Let us denote by $\mathbb{G}_{(L)}$ (resp. $\mathbb{G}_\infty$) the largest ($\Vert \cdot\Vert_1$-)connected component of $\mathbb{Y}$ contained in $[-L,L]^d\cap\ZZ^d$ (resp. the a.s. unique infinite component of $\mathbb{Y}$). We then define $\mathcal{G}_{(L)}=\mathcal{G}_{(L)}(\widehat{\xi})$ (resp. $\mathcal{G}_{\infty}=\mathcal{G}_\infty(\widehat{\xi})$) as the set of points of $\xi$ whose Voronoi cell intersects a $K-$box with index in $\mathbb{G}_{(L)}$ (resp. $\mathbb{G}_\infty$). The points of $\mathcal{G}_\infty$ are called \emph{good points}. Let us note that $\mathcal{G}_{(L)}$ and $\mathcal{G}_\infty$ are connected in the Delaunay triangulation of $\xi$.

We claim that good points have their degrees uniformly bounded. More precisely:
\begin{lemm}\label{LemmBoundedDegree}
There exists $D=D(d, K, \alpha)$ such that for every $x\in \xi$ with Voronoi cell intersecting an $\alpha$-good box:
\[\deg_{\DT(\xi)}(x)\leq D.\] 
\end{lemm}
\begin{rem}
As it appears in the following proof, $D$ is also an upper bound for the maximal number of Voronoi cells which intersect a good box. It will be used to bound these two quantities throughout the paper.
\end{rem}
\begin{dem}
Let $x\in \xi$ be such that its Voronoi cell intersects an $\alpha$-good box $B_\z$. By definition, each box $B_\z'$ with $\Vert \z'-\z\Vert_\infty\leq 1$ is $\alpha$-nice. Since any sub-box of side $s$ that belongs to a nice box contains at least one point of $\xi$, the points in nice boxes are whithin a distance at most $\sqrt{d}s$ from the nucleus of their Voronoi cell. In particular, $x$ is whithin a distance at most $\sqrt{d}s$ from $B_\z$. Similarly, the nuclei of the Voronoi cells which share a face with $\Vor_\xi(x)$ are in $B_\z+B_2(0, 3\sqrt{d}s)\subset \overline{B_\z}:=\bigcup_{\z':\Vert \z'-\z\Vert_\infty\leq 1}B_{\z'}$. Hence, $\deg_{\DT(\xi)}(x)$ is generously bounded by $\#\big(\xi\cap \overline{B_\z}\big)\leq D:=\alpha (3K)^d$.
\end{dem}
\subsubsection{The ergodicity assumption}\label{ErgAss} As in \cite{CFP}, we have an ergodicity assumption but we make more explicit the use of the coupling. We will use it to adapt the method developed in \cite[\S 4]{BergerBiskup} for proving the sublinearity of the corrector along coordinate directions. 

\begin{itemize}
\item[{\bf (Er)}] For each $(K,p)$ such that {\bf (SD)} holds and each $e$, $\P=\P^{K,e}$ is ergodic with respect to the transformation:
\[\tau=\tau^{K,e} : (\xi, (y_\z)_{\z\in\ZZ^d})\longrightarrow (\tau_{Ke}\xi,(y_{\z+e})_{\z\in\ZZ^d}).\]
\end{itemize}

Note that $\tau$ is invertible and $\P$ is invariant with respect to $\tau$ due to the stationarity of the point process.

\subsubsection{Polynomial moments} We also make the following assumption:
\begin{itemize}
\item[{\bf (PM)}]
the number of points in a unit cube admits a polynomial moment of order 2 under $\P$ and $\deg_{\DT(\xi^0)}(0)$ and $\max_{x\sim 0\text{ in }\DT(\xi^0)}\Vert x\Vert$ admit respectively a moment of order 2 and 4 under the Palm distribution $\P_0$ associated with the point process.
\end{itemize}
\subsubsection{The case of point processes with a finite range of dependence}
We verify that, when the point process has a finite range of dependence, assumption  {\bf (Er)} is always satisfied, and assumptions {\bf (SD)} and {\bf (PM)} are implied by assumptions {\bf (V)}, {\bf (D)}, {\bf (V')} and {\bf (EM)} which are described below.

Let us check that {\bf (Er)} is satisfied if the point process has a finite range of dependence. Let $\A$ be a measurable set such that $\tau\A=\A$. Fix $\varepsilon>0$ and $\B$ with $\P(\A\Delta\B)\leq\varepsilon$ which depends only on $\xi $ restricted to a compact subset of $\RR^d$ and on finitely many $y_\z$s. It is thus possible to find $n$ such that $\B$ and $\tau^n\B$ are independent. Then, $\P[\B\cap\tau^n\B]=\P[\B]^2$ by invariance of $\P$ w.r.t. $\tau$. Hence, 
\begin{align*}
\big\vert \P[\A]-\P[\A]^2\big\vert &\leq \big\vert \P[\A]-\P[\B\cap \tau^n\B]\big\vert+\big\vert \P[\B\cap \tau^n\B]-\P[\A]^2\big\vert\\
 &\leq \P[\A\Delta (\B\cap \tau^n\B)]+\big\vert \P[\B]^2-\P[\A]^2\big\vert\\
 &\leq \P[(\A\Delta \B)\cup (\A\Delta \tau^n\B)]+\big( \P[\B]+\P[\A]\big)\big\vert \P[\B]-\P[\A]\big\vert\\
 &\leq \P[\A\Delta \B]+\P[\A\Delta \tau^n\B]+2\varepsilon\\
 &\leq \P[\tau^n\A\Delta \tau^n\B]+3\varepsilon\\
 &= \P[\A\Delta \B]+3\varepsilon\leq 4\varepsilon
\end{align*}
Since $\varepsilon$ is arbitrary, this implies that $\P[\A]$ is 0 or 1.
\bigskip

Since deciding if a box is good depends only on the behavior of the point process $\xi$ in a neighborhood of the box, the process of good boxes $\mathbb{X}$ is a Bernoulli process on $\ZZ^d$ with a finite range of dependence when the point process has itself a finite range of dependence. In this case, if $\P[X_\z=1]$ is as close to 1 as we wish for $s$ and $\alpha$ large enough, \cite[Theorem 0.0]{Liggett97} ensures that {\bf (SD)} holds. One can easily bound $\P[X_\z=1]$ from below when the point process satisfies {\bf (V)} and
\begin{itemize}
\item[{\bf (D)}] \label{condproc2} 
there exist positive constants $c_2,c_3$ such that for $L$ large enough:
\begin{equation*}
\P\big[\#\big([0,L]^d\cap\xi\big)\geq c_2L^d\big]\leq e^{-c_3L^d}.
\end{equation*}
\end{itemize}

In Lemma \ref{MomentsExpo}, we prove that $\deg_{\DT(\xi^0)}(0)$ and $\max_{x\sim 0\text{ in }\DT(\xi^0)}\Vert x\Vert$ admit exponential moments when the point process has a finite range of dependence and satisfies: 
\begin{itemize}
\item[{\bf (V')}]
there exists a positive constant $c_4$ such that for $L$ large enough:
\begin{equation*}\label{MomentsExpoVoid}
\P_0\Big[\#\big(\xiz\cap C_L\big)=0\Big]\leq e^{-c_4L^d},\quad \forall C_L\mbox{ cube of side }L, 
\end{equation*}
\end{itemize}
and
\begin{itemize}
\item[{\bf (EM)}] 
there exist a positive constant $c_5$ and a positive function $f(\rho)$ which goes to 0 with $\rho$ such that for $L$ large enough:
\begin{equation*}\label{MomentsExpoPalm}
\E_0\Big[e^{\rho\#(\xiz \cap [-L,L]^d)}\Big]\leq c_5e^{f(\rho)L^d}.
\end{equation*}
\end{itemize}

Let us finally note that these assumptions are in particular satisfied by homogeneous Poisson point processes, Mat\'ern cluster processes and type I or II Mat\'ern hardcore processes. Indeed, these point processes have finite range of dependence and it is quite classical to check assumptions {\bf (V)}, {\bf (D)}, {\bf (V')} and {\bf (EM)} for these processes (see \cite[Appendix]{AIP}).
\subsection{Outline of the paper} As announced at the beginning of the introduction, the crux is to prove Theorem \ref{TheoSousLine} and we follow the approach of \cite{BP}. The main steps of the proof are stated explicitly in Theorem 2.4 of that paper. We prove the existence of the corrector in Section \ref{SectCorr}. Next, we verify that the corrector grows at most polynomially in Section \ref{SectPolGrowth} and at most linearly in each coordinate direction in Section \ref{DirSubli}. The sublinearity on average is treated in Section \ref{AverSubli} while diffusive bounds for a related random walk are proved in Sections \ref{HeatKernelBound} and \ref{ExpectedDistanceBound}. In Section \ref{ASRestSubli}, we prove the a.s. sublinearity of the corrector in the set of `good points'. The proofs  of Theorems \ref{QIP} and \ref{TheoSousLine} are finally completed in Section \ref{SubliProof}.
\section{Construction of the corrector and harmonic deformation}\label{SectCorr}
Let us define the measure $\mu$ on $\N_0\times \RR^d$ by:
\[\int u\d\mu:=\E_0\Big[\sum_{x\in\xi^0}c^{\xi^0}_{0,x}u(\xi^0,x)\Big],\]
where $c^{\xi^0}_{0,x}=\mathbf{1}_{0\sim x\text{ in }\DT(\xi^0)}$.
This measure has total mass $\E_0[\deg_{\DT(\xi^0)}(0)]$ which is finite thanks to assumption {\bf (PM)}.  We denote by $(\cdot,\cdot)_\mu$ the scalar product in $L^2(\mu)$. 
\subsection{Weyl decomposition of $L^2(\mu)$}
As in \cite{MP}, \cite{BergerBiskup} or \cite{CFP}, we work with the orthogonal decomposition of $L^2(\mu)$ in the subspaces of square integrable \emph{potential} and \emph{solenoidal} fields. This decomposition is quite standard (see \emph{e.g.} \cite[Chap. 9]{Lyons}) and generally called \emph{Weyl decomposition}.  Let us denote by $(\tau_x)_{x\in \RR^d}$ the group of translations in $\RR^d$ which acts naturally on $\N_0$ as follows: $\tau_x\xi^0=\sum_{y\in\xi^0}\delta_{y-x}$.

\begin{defi} For $\psi : \N_0\rightarrow \RR$, the \emph{gradient field} $\nabla \psi:\,\N_0\times \RR^d\longrightarrow \RR$ is defined for $x\in\xi^0$ by:
\[\nabla \psi (\xi^0,x):=\psi (\tau_x\xi^0)-\psi(\xi^0)\]
and by 0 if $x\not\in\xi^0$.
\end{defi}
Note that gradients of measurable bounded functions on $\N_0$ are elements of $L^2(\mu)$ thanks to assumption {\bf (PM)}.
\begin{defi}The space $L^2_\text{pot}(\mu)$ of \emph{potential} fields is defined as the closure of the subspace of gradients of measurable bounded functions on $\N_0$. Its orthogonal complement is the set of \emph{solenoidal} (or \emph{divergence-free}) fields and is denoted by $L^2_\text{sol}(\mu)$.
\end{defi}

Let us recall some additional definitions:
\begin{defi}\label{DefiAntiEtc}
A function $u:\,\N_0\times\RR^d\longrightarrow \RR$ is called:
\begin{enumerate}
\item \emph{antisymmetric} if
\[u(\xi^0,x)=-u(\tau_x\xi^0,-x), \qquad \forall \xi^0\in\N_0,\,\forall x\in\xi^0;\]
\item \emph{shift-covariant} if
\[u(\xi^0,x)=u(\xi^0,y)+u(\tau_y\xi^0,x-y), \qquad \forall \xi^0\in\N_0,\,\forall x,y\in\xi^0;\]
\item \emph{curl-free} if it satisfies the following \emph{co-cycle} relation: for any $\xi^0\in\N_0$, any $n\in\NN^*$, and any collection of points $x_0, \dots,x_n\in\xi^0$ with $x_0=x_n$, one has
\begin{equation}\label{EqDefCurlFree}
\sum_{i=1}^{n-1}u(\tau_{x_i}\xi^0,x_{i+1}-x_i)=0.
\end{equation}
\end{enumerate}
\end{defi}
A function of $L^2(\mu)$ is called \emph{antisymmetric} (resp. \emph{shift-covariant}, \emph{curl-free}) if it is antisymmetric (resp. shift-covariant, curl-free) for $\P_0$-a.a. $\xi^0$. In each case, it admits a representative wich satisfies the corresponding property everywhere.
By taking $x=y=0$ in the definition, one can see that any shift-covariant function $u$ must satisfy $u(\xi^0,0)=0$ for any $\xi^0\in\N_0$. Next proposition lists simple but useful links between definitions above:
\begin{prop}
Let $u\in L^2(\mu)$.
\begin{enumerate}
\item If $u\in L^2_\text{pot}(\mu)$, then it is curl-free.  
\item If $u$ is curl-free, then it is also antisymmetric and shift-covariant.
\end{enumerate}
\end{prop}
\begin{dem}$\,$
\begin{enumerate}
\item Gradients fields are clearly curl-free. The general case is obtained by a standard approximation argument.
\item The case $n=2$, $x_0=x_2=0$ and $x_1=x$ in (\ref{EqDefCurlFree}) gives the antisymmetry of $u$. By (\ref{EqDefCurlFree}) with $n=3$, $x_0=x_3=0$, $x_1=y$ and $x_2=x$, one has:
\[u(\xi^0,y)+u(\tau_y\xi^0, x-y)+u(\tau_x\xi^0, -x)=0.\] 
The shift-covariance of $u$ then follows from its antisymmetry.
\end{enumerate}
\end{dem}

We can now define the \emph{divergence} of an integrable field and derive an integration by parts formula.

\begin{defi}
The \emph{divergence} of $u\in L^1(\mu)$ is defined by:
\[\div u(\xi^0):=\sum_{x\in\xi^0}c^{\xi^0}_{0,x}u(\xi^0,x),\qquad \xi^0\in\N_0.\]
\end{defi}
Triangle inequality clearly implies that divergences of $L^1(\mu)$ functions are in $L^1(\P_0)$. Moreover, if $u\in L^1(\mu)$ is a positive function, we have the equality:
\begin{equation}\label{EqnormesL1}
\Vert u \Vert_{L^1(\mu)}=\Vert \div u \Vert_{L^1(\P_0)}.
\end{equation}
We derive the following integration by parts formula:
\begin{lemm} Let $\psi$ be a bounded measurable function on $\N_0$ and let $u\in L^2(\mu)$ be an antisymmetric field. It holds that:
\begin{equation}\label{IPP}
(u,\nabla\psi )_\mu=-2\E_0[\psi\div u].
\end{equation}
\end{lemm}
\begin{dem} Observe that:
\[c^{\xi^0}_{0,x}=c^{\tau_x\xi^0}_{0,-x}.\]
Due to the antisymmetry of $u$, one has:
\begin{align*}
(u,\nabla\psi )_\mu &=\E_0\Bigg[\sum_{x\in\xi^0}c^{\xi^0}_{0,x}u(\xi^0,x)\psi (\tau_x\xi^0)\Bigg]-\E_0\Bigg[\sum_{x\in\xi^0}c^{\xi^0}_{0,x}u(\xi^0,x)\psi (\xi^0)\Bigg]\\
&=-\E_0\Bigg[\sum_{x\in\xi^0}c^{\tau_x\xi^0}_{0,-x}u(\tau_x\xi^0,-x)\psi (\tau_x\xi^0)\Bigg]-\E_0[\psi\div u]. 
\end{align*}
By Neveu exchange formula (see \cite[Theorem 3.4.5]{SW} in the special case where $X=Y$), one has for any integrable function $f$ on $\N_0\times \RR^d$:
\[\int_{\N_0}\sum_{x\in\xi^0}f(\tau_x\xi^0,x)\P_0(\d\xi^0)=\int_{\N_0}\sum_{x\in\xi^0}f(\xi^0,-x)\P_0(\d\xi^0).\] 
The conclusion is then obtained by applying the identity above with: $$f(\xi^0,x):=c^{\xi^0}_{0,-x}u(\xi^0,-x)\psi(\xi^0).$$
\end{dem}

This lemma implies the immediate following corollary.
\begin{coro}\label{CorDiv1}An antisymmetric field $u\in L^2(\mu)$ is solenoidal if and only if $\div u(\xi^0)=0$ for $\P_0$-a.e. $\xi^0$.
\end{coro}
\subsection{Construction of the corrector}
In \cite{FGG}, the authors relied on an Harness-type process to obtain the existence harmonic deformations of Delaunay triangulations which corresponds to the existence of the corrector. Here, we recall how the decomposition of $L^2(\mu)$ allows us to derive the existence of the corrector by following the construction of \cite{MP,CFP}.

For $i=1,\dots, d$, $\xi^0\in\N_0$ and $x\in\RR^d$, set $u_i(\xi^0, x):=x_i$ the $i^\text{th}$ coordinate of $x$. Note that $u_i$ is clearly antisymmetric in the sense of Definition \ref{DefiAntiEtc} and that $u_i\in L^2 (\mu)$. Actually, by the Cauchy-Schwarz inequality and assumption {\bf (PM)}, one obtains:
\begin{align*}\int \vert u_i\vert^2\d\mu=\E_0\Bigg[\sum_{x\in\xi^0}c^{\xi^0}_{0,x}\vert x_i\vert^2\Bigg]
&\leq\E_0\Big[\big(\max_{x\sim 0\tiny\text{ in }\DT (\xiz)}\Vert x\Vert\big)^2\deg_{\DT(\xiz)}(0)\Big]\\
&\leq\E_0\Big[\big(\max_{x\sim 0\tiny\text{ in }\DT (\xiz)}\Vert x\Vert\big)^4\Big]^\frac{1}{2}\E_0\Big[\big(\deg_{\DT(\xiz)}(0)\big)^2\Big]^\frac{1}{2}<\infty.
\end{align*}
Consider now the orthogonal decomposition of the form $u_i=\chi_i+\varphi_i$ with $\chi_i\in L^2_\text{pot}(\mu)$ and $\varphi_i\in L^2_\text{sol}(\mu)$. Since $\chi_i\in L^2_\text{pot}(\mu)$, it is antisymmetric and $\varphi_i$ also (as a difference of antisymmetric functions). Hence, it follows from Corollary \ref{CorDiv1} that $\varphi_i$ is harmonic at 0.

The \emph{corrector field} is the vector-valued function on $\N_0\times \RR^d$ defined by	$\chi=(\chi_1,\dots,\chi_d)$. It admits a shift-covariant representative and its norm $\Vert\chi\Vert$ is in $L^2(\mu)$. We denote by $\varphi=(\varphi_1,\dots,\varphi_d)$ the harmonic function $(\xi^0,x)\longrightarrow x-\chi(\xi^0,x)$. From the harmonicity of $\varphi(\xi^0,\cdot)$ at 0 for $\P_0$-a.a. $\xi^0$ and \cite[Lemma B.2]{CFP}, one deduces that, for $\P$-a.a. $\xi$, for any $x\in\xi$:
\begin{equation}\label{EqHarmo}
\sum_{y\in\xi}c^{\xi}_{x,y}\Vert \varphi(\tau_x\xi,y-x)\Vert<\infty\qquad\mbox{and}\qquad\sum_{y\in\xi}c^{\xi}_{x,y}\varphi(\tau_x\xi,y-x)=0.
\end{equation}
Let us define $$M^\xi_n:=\varphi(\tau_{X_0}\xi,X_n-X_0)=\sum_{i=1}^n\varphi(\tau_{X_i}\xi,X_{i+1}-X_i).$$ It follows from (\ref{EqHarmo}) that $(M^\xi_n)_{n\in\NN}$ is a martingale under $P^\xi_x$.
\section{Polynomial growth}\label{SectPolGrowth}
Let us define:
\begin{equation}\label{EqCroissPol1}
\mathcal{R}_n=\mathcal{R}_n(\xi):=\max_{\tiny\begin{array}{c}x,y\in\xi\\ 
\Vor_{\xi}(x)\cap[-n,n]^d\neq\emptyset\\
\Vor_{\xi}(y)\cap[-n,n]^d\neq\emptyset\end{array}}\big\Vert \chi(\tau_x\xi,y-x)\big\Vert.
\end{equation}
\begin{prop}\label{LemmCroissPol}
For every $\beta >d+1$, one has:
\begin{equation*}\label{EqCroissPol2}
\frac{\mathcal{R}_n}{n^\beta} \xrightarrow[n\rightarrow \infty]{\,\P-a.s.\,}0.
\end{equation*}
\end{prop}
\begin{dem} For $n$ fixed, let us cover $[-2n,2n]^d$ with disjoint boxes of side $\log n$ and denote by $\mathcal{A}_n$ the event `each of these boxes contains at least one point of $\xi$'. Note that, thanks to assumption {\bf (V)},  $\P[\mathcal{A}_n^c]=O(n^{-2})$. Let $d+1<{\beta'}<\beta$. Since $$\P[\mathcal{R}_n\geq n^{\beta'}]\leq \P[\mathcal{R}_n\mathbf{1}_{\mathcal{A}_n}\geq n^{\beta'}]+\P [\mathcal{A}_n^c],$$ we only need to show that $\sum_n\P[\mathcal{R}_n\mathbf{1}_{\mathcal{A}_n}\geq n^{\beta'}]<\infty$. The result then follows using the Borel-Cantelli lemma.

Let $\x,\y\in\xi$ with $\Vor_{\xi}(\x)\cap[-n,n]^d\neq\emptyset$ and $\Vor_{\xi}(\y)\cap[-n,n]^d\neq\emptyset$ be such that $\mathcal{R}_n=\Vert \chi (\tau_\x\xi, \y-\x)\Vert$. Consider the simple Delaunay-path $(x_0, \dots, x_m)$ from $x_0=\x$ to $x_m=\y$ obtained by connecting the nuclei of successive Voronoi cells which intersect the line segment $[\x,\y]$. Observe that, on $\mathcal{A}_n$, any point of $[-2n,2n]^d$ is within a distance at most $\sqrt{d}\log n$ from the nucleus of its Voronoi cell. In particular, $(x_0, \dots, x_m)$ is contained in $[-2n,2n]^d$. 

Recall that the corrector $\chi$ is shift-covariant. Hence, for $i=0,\dots, m-1$, one has:
\[\chi (\tau_{x_i}\xi,x_{i+1}-x_i)=\chi (\tau_\x\xi,x_{i+1}-\x)-\chi (\tau_\x\xi,x_{i}-\x),\]
and using that $\chi (\tau_\x\xi,0)=0$:
\[\chi (\tau_\x\xi,\y-\x)=\sum_{i=1}^{m-1}\chi (\tau_{x_i}\xi,x_{i+1}-x_i).\]

We deduce that on $\mathcal{A}_n$:
\begin{align*}
\mathcal{R}_n&\leq \sum_{i=1}^{m-1}\Vert \chi (\tau_{x_i}\xi,x_{i+1}-x_i)\Vert\\
&\leq \sum_{x\in\xi\cap [-2n,2n]^d}\sum_{y\in\xi}c_{x,y}^\xi \Vert \chi (\tau_{x}\xi,y-x)\Vert\\
&= \sum_{x\in\xi\cap [-2n,2n]^d}\div \Vert \chi \Vert  (\tau_{x}\xi).
\end{align*}
Together with Markov inequality and Campbell formula, the inequality above leads to:
\begin{align*}
\P[\mathcal{R}_n\mathbf{1}_{\mathcal{A}_n}\geq n^{\beta'}]&\leq \frac{\E[\mathcal{R}_n\mathbf{1}_{\mathcal{A}_n}]}{n^{\beta'}}\\
&\leq n^{-\beta'}\int_\N\sum_{x\in\xi\cap [-2n,2n]^d}\div \Vert \chi \Vert  (\tau_{x}\xi)\P(\d\xi)\\
&\leq cn^{d-\beta'}\big\Vert\div \Vert \chi \Vert \big\Vert_{L^1(\P_0)}
\end{align*}
Since ${\beta'}>d+1$ and $\div \Vert \chi \Vert \in L^2(\P_0)$, this completes the proof. 
\end{dem}
\section{Sublinearity along coordinate directions in $\mathcal{G}_\infty(\widehat{\xi})$}\label{DirSubli}
In this section,  we adapt the arguments of \cite[\S 7.2]{CFP} which consist in an adaptation of the `lattice method' developed in \cite{BergerBiskup,BP}.

Given a unit vector $e$ in the direction of one of the coordinate axes of $\RR^d$ and $\widehat{\xi}\in\widehat{\N}$, let us define:
\[n_0(\widehat{\xi})=n_0^e(\widehat{\xi}):=0\qquad\mbox{and}\qquad n_{i+1}(\widehat{\xi})=n_{i+1}^e(\widehat{\xi}):=\min\{j>n_i(\widehat{\xi})\,: je\in \mathbb{G}_\infty\}.\]

Recall the definition of $\tau$ from assumption {\bf{(Er)}} and consider the shift $\tau_*$ induced on $\N^*:=\{\widehat{\xi}\in\widehat{\N}:0\in\mathbb{G}_\infty\}$ from $\tau$, that is
$\tau_*:\widehat{\xi}\longrightarrow \tau^{n_1(\widehat{\xi})}\widehat{\xi}$.  
Thanks to assumption {\bf{(Er)}} and the fact that $\P\left[0\in\mathbb{G}_\infty\right]>0$, standard arguments (see {\it e.g.} \cite[Lemma 7.3]{CFP} or \cite[Theorem 3.2]{BergerBiskup}) lead to:
\begin{lemm}\label{LemmErgo*}The probability measure $\P\left[\cdot\vert 0\in\mathbb{G}_\infty\right]$ is stationary and ergodic w.r.t. $\tau_*$.
\end{lemm}

Next, for $\widehat{\xi}\in\mathcal{N}^*$, we write $w_i$ for the point of $\xi$ whose Voronoi cell contains the center of the box $B_{n_i(\widehat{\xi})Ke}$.

\begin{lemm}\label{LemmIntChi*}It holds that:
 \[\E\left[\left\Vert\chi(\tau_{w_0}\xi, w_1-w_0)\right\Vert\big\vert 0\in \mathbb{G}_\infty\right]<\infty\qquad \mbox{and}\qquad \E\left[\chi(\tau_{w_0\xi}, w_1-w_0)\big\vert 0\in \mathbb{G}_\infty\right]=0.\]
\end{lemm}
\begin{dem} For $\widehat{\xi}\in\N^*$, let $d(\widehat{\xi})\geq n_1(\widehat{\xi})$ denote the chemical distance between 0 and $n_1(\widehat{\xi})e$ in the infinite cluster $\mathbb{G}_\infty$. On the event $\{d(\widehat{\xi})=j\}$, there exists a path $\z_0=0,\z_1,\dots,\z_j=n_1(\widehat{\xi})e$ in $\mathbb{G}_\infty$. For $i=0,\dots, j-1$, thanks to the definition of the good boxes and assumption {\bf (SD)}, the nuclei of the Voronoi cells intersecting the line segment $[K\z_i,K\z_{i+1}]$ are within a distance at most $\sqrt{d}s$ from this line segment.  By connecting the successive nuclei of the Voronoi cells intersecting the broken line $[K\z_0,K\z_1,\dots ,K\z_j]$, we obtain a simple path $w_0=x_0,x_1,\dots,x_m=w_1$ between $w_0$ and $w_1$ in $\mathcal{G}_\infty(\widehat{\xi})$ which is contained in $[-K(j+\frac{1}{2}),K(j+\frac{1}{2})]^d$ and has length $m\leq (j+1)D$. Thanks to the shift-covariance of $\chi$, as in the proof of Proposition \ref{LemmCroissPol}, one obtains:
\begin{align*}
 \Vert\chi(\tau_{w_0}\xi, w_1-w_0)\Vert
 &\leq\sum_{i=1}^{m-1}\Vert\chi(\tau_{x_i}\xi, x_{i+1}-x_i)\Vert\\
 &\leq\sum_{\tiny\begin{array}{c}x\in\mathcal{G}_\infty(\widehat{\xi}):\\\Vert x\Vert_\infty \leq K(j+\frac{1}{2})\end{array}}\sum_{y\underset{\DT (\xi)}{\sim} x}\Vert\chi(\tau_{x}\xi, y-x)\Vert.
\end{align*}
Together with the Cauchy-Schwarz inequality, this leads to:
\begin{align}\label{EqIntInduCor1}
\E[\Vert&\chi(\tau_{w_0}\xi, w_1-w_0)\Vert\big\vert 0\in\mathbb{G}_\infty]\nonumber\\&=\sum_{j=1}^\infty \E\left[\Vert\chi(\tau_{w_0}\xi, w_1-w_0)\Vert\mathbf{1}_{d(\widehat{\xi})=j}\big\vert 0\in\mathbb{G}_\infty\right]\nonumber\\
&\leq\sum_{j=1}^\infty \E\Bigg[\sum_{\tiny\begin{array}{c}x\in\mathcal{G}_\infty(\widehat{\xi}):\\\Vert x\Vert_\infty \leq K(j+\frac{1}{2})\end{array}}\sum_{y\underset{\DT (\xi)}{\sim} x}\Vert\chi(\tau_{x}\xi, y-x)\Vert\mathbf{1}_{d(\widehat{\xi})=j}\Big\vert 0\in\mathbb{G}_\infty\Bigg]\nonumber\\
%
&\leq\sum_{j=1}^\infty \Bigg\{\E\Bigg[\sum_{\tiny\begin{array}{c}x\in\mathcal{G}_\infty(\widehat{\xi}):\\\Vert x\Vert_\infty \leq K(j+\frac{1}{2})\end{array}}\sum_{y\underset{\DT (\xi)}{\sim} x}\Vert\chi(\tau_{x}\xi, y-x)\Vert^2\Big\vert 0\in\mathbb{G}_\infty\Bigg]^\frac{1}{2}\nonumber\\&\qquad\qquad\qquad\times\E\Bigg[\Bigg(\sum_{\tiny\begin{array}{c}x\in\mathcal{G}_\infty(\widehat{\xi}):\\\Vert x\Vert_\infty \leq K(j+\frac{1}{2})\end{array}}\deg_{\DT(\xi)}(x)\Bigg)\mathbf{1}_{d(\widehat{\xi})=j}\Big\vert 0\in\mathbb{G}_\infty\Bigg]^\frac{1}{2}\Bigg\}\nonumber\\
&\leq\frac{1}{\P\big[0\in\mathbb{G}_\infty\big]^\frac{3}{4}}
\sum_{j=1}^\infty \Bigg\{\E\Bigg[\sum_{\tiny\begin{array}{c}x\in\xi:\\\Vert x\Vert_\infty \leq K(j+\frac{1}{2})\end{array}}\div\Vert\chi\Vert^2(\tau_x\xi)\Bigg]^\frac{1}{2}\nonumber
\\&\qquad\qquad\qquad\times
\E\Bigg[\Bigg(\sum_{\tiny\begin{array}{c}x\in\mathcal{G}_\infty(\widehat{\xi}):\\\Vert x\Vert_\infty \leq K(j+\frac{1}{2})\end{array}}\deg_{\DT(\xi)}(x)\Bigg)^2\Bigg]^\frac{1}{4}
\P\big[d(\widehat{\xi})=j\big\vert 0\in\mathbb{G}_\infty\big]^\frac{1}{4}\Bigg\}.
\end{align}

It follows from Campbell formula and formula (\ref{EqnormesL1}) that $$\E\left[\sum_{\tiny x\in\xi:\Vert x\Vert_\infty \leq K(j+\frac{1}{2})}\div\Vert\chi\Vert^2(\tau_x\xi)\right]\leq cj^d\big\Vert \div\Vert\chi\Vert^2 \big\Vert_{L^1(\P_0)}=cj^d\big\Vert \Vert\chi\Vert \big\Vert^2_{L^2(\mu)}.$$ Since points of $\mathcal{G}_\infty(\widehat{\xi})$ have degrees bounded by $D$ (see Lemma \ref{LemmBoundedDegree}) and $\#(\xi\cap [0,1]^d)$ admits a moment of order 2, one has:
\[\E\Bigg[\Bigg(\sum_{\tiny\begin{array}{c}x\in\mathcal{G}_\infty(\widehat{\xi}):\\\Vert x\Vert_\infty \leq K(j+\frac{1}{2})\end{array}}\deg_{\DT(\xi)}(x)\Bigg)^2\Bigg]\leq D^2\E\Big[\# \Big(\xi\cap\Big[-K\Big(j+\frac{1}{2}\Big),K\Big(j+\frac{1}{2}\Big)\Big]^d\Big)^2\Big]\leq cj^d.\]

Thanks to \cite[Lemma 4.4]{BergerBiskup}, we know that $\P\big[d(\widehat{\xi})=j\big\vert 0\in\mathbb{G}_\infty\big]$ has an exponential decay. Hence, collecting bounds, we obtain that the sum in the r.h.s. of (\ref{EqIntInduCor1}) is finite.

Since $\chi\in L^2_\text{pot}(\mu)$, it is the $L^2$-limit of gradients of bounded measurable functions $(g_n)_{n\in\NN}$ defined on $\N_0$. Let us define $\chi_n:=\nabla g_n$. By the same arguments as above, one obtains that:
\[\E\big[\Vert \chi(\tau_{w_0}\xi, w_1-w_0)-\chi_n(\tau_{w_0}\xi, w_1-w_0)\Vert\big\vert 0\in\mathbb{G}_\infty\big]\leq c\Vert \chi-\chi_n\Vert^2_{L^2(\mu)}\xrightarrow[\, n\rightarrow\infty\,]{}0.\]

Note that, for all $i$, $w_i$ is a deterministic function of $\widehat{\xi}$ and that $w_1(\widehat{\xi})=w_0(\tau_*\widehat{\xi})+n_1(\widehat{\xi})Ke$.
The conclusion follows since, due to the stationarity of $\P[\cdot\vert 0\in\mathbb{G}_\infty]$ with respect to $\tau_*$ applied to the function $\left(g_n\circ \tau_{w_0(\cdot)}\right)(\cdot)$, $\E\big[\chi_n(\tau_{w_0}\xi, w_1-w_0)\big\vert 0\in\mathbb{G}_\infty\big]=0$. 
\end{dem}

Combining Lemmas \ref{LemmErgo*} and \ref{LemmIntChi*}, one obtains the sublinearity along the direction $e$ in $\mathcal{G}_\infty(\widehat{\xi})$.
\begin{prop}\label{LemmeSousLineariteDir}
For $\P[\cdot\vert 0\in\mathbb{G}_{\infty}]-a.a.\,\widehat{\xi}$:
\[\lim_{k \rightarrow\infty}\frac{\chi(\tau_{w_0}\xi, w_k-w_0)}{k}=0,\]
and
\[\lim_{k \rightarrow\infty}\max_{\tiny\begin{array}{c}
x_0\in \xi :\\\Vor_\xi(x_0)\cap B^K_0\neq\emptyset\end{array}}\frac{\chi(\tau_{x_0}\xi, w_k-x_0)}{k}=0.\]
\end{prop}
\begin{dem} Thanks to the shift-covariance of the corrector, one has:
\begin{align}\label{EqDirSub2}
\frac{\chi(\tau_{w_0}\xi, w_k-w_0)}{k}&=\frac{1}{k}\sum_{j=0}^{k-1}\big(\chi(\tau_{w_0}\xi, w_{j+1}-w_0)-\chi(\tau_{w_0}\xi, w_j-w_0)\big)\nonumber\\
&=\frac{1}{k}\sum_{j=0}^{k-1}\chi(\tau_{w_j}\xi,w_{j+1}-w_{j}).
\end{align}
Observe that $w_j(\widehat{\xi})=w_0(\tau_*^j\widehat{\xi})+n_j(\widehat{\xi})Ke$, $w_{j+1}(\widehat{\xi})=w_1(\tau_*^j\widehat{\xi})+n_1(\tau_*^j\widehat{\xi})Ke$ and recall that, by Lemma \ref{LemmErgo*}, $\P[\cdot\vert 0\in\mathbb{G}_\infty]$ is ergodic with respect to $\tau_*$. Thanks to Birkhoff's theorem, the last expression in (\ref{EqDirSub2}) converges to $\E[\chi(\tau_{w_0}\xi,w_{1}-w_0)\vert 0\in\mathbb{G}_\infty]$ which is 0 by Lemma \ref{LemmIntChi*}.

The second part of the lemma then follows from equality :
\[\chi(\tau_{x_0}\xi, w_k-x_0)=\chi(\tau_{w_0}\xi, w_k-w_0)+\chi (\tau_{x_0}\xi, w_0-x_0)\]
which is due to the shift-covariance of the corrector.
\end{dem}
\section{Sublinearity on average in $\mathcal{G}_\infty(\widehat{\xi})$}\label{AverSubli}
We derive the sublinearity on average of the corrector in $\mathcal{G}_\infty(\widehat{\xi})$ from Lemma \ref{LemmeSousLineariteDir}. Our approach is close in spirit to \cite[\S 7.3]{CFP}.
\begin{prop}\label{PropSousLineariteMoy}
For every $\varepsilon_0>0$, for $\P[\cdot \vert 0\in\mathbb{G}_\infty]-a.a.~\widehat{\xi}\in\widehat{\N}$:
\begin{equation}
\lim_{L\rightarrow \infty}\max_{\tiny\begin{array}{c}x_0\in \xi :\\\Vor_\xi(x_0)\cap B^K_0\neq\emptyset\end{array}}\frac{1}{\#\Lambda_L}\sum_{\tiny\begin{array}{c}x\in \xi :\\
\exists \z\in\mathbb{G}_\infty\cap \Lambda_L\\
\Vor_\xi(x)\cap B^K_\z\neq\emptyset\end{array}}\mathbf{1}_{\Vert\chi (\tau_{x_0}\xi,x_0-x)\Vert\geq \varepsilon_0L}=0,
\end{equation}
where $\Lambda_L:=\ZZ^d\cap [-L,L]^d$.
\end{prop}

Let us describe roughly the method which is an alternative to the one of \cite[\S 5.2]{BergerBiskup} and relies on multiscale arguments. The first idea is to extend the directional sublinearity result of Proposition \ref{LemmeSousLineariteDir} dimension by dimension. For $\nu\in\{1,\dots,d\}$, we denote by $\Lambda^\nu_L$ the set: 
\[\Lambda^\nu_L:=\left\{\z=(\z_1,\dots,\z_d)\in\ZZ^d:~\mbox{for } 1\leq i\leq\nu, \vert\z_i\vert\leq L, \mbox{ and for }\nu+1\leq i\leq d, \z_i=0\right\}.\]
Assume that we have a `good' (sublinear) control of $\chi (\tau_{x_0}\xi, x-x_0)$ for $x_0$ whose Voronoi cell intersects $B^K_0$ and $x$ whose Voronoi cell intersects $B^K_\z$ for some $\z\in \Lambda_L^{\nu}\cap \mathbb{G}_\infty$. Then, using Proposition \ref{LemmeSousLineariteDir}, one obtains a sublinear control on $\chi (\tau_{x}\xi, x'-x)$ for $x$ whose Voronoi cell intersects $B^K_\z$ and $x'$ whose Voronoi cell intersects $B^K_{\z'}$ for any $\z'\in \Lambda_L^{\nu+1}\cap \mathbb{G}_\infty$ which differs from $\z$ only on the $(\nu+1)$-th coordinate. By the shift-covariance, this gives a control on $\chi (\tau_{x_0}\xi, x'-x_0)$. As noticed in \cite{BergerBiskup}, we can not deduce directly Proposition \ref{PropSousLineariteMoy} from this argument because $\mathbb{G}_\infty$ covers only a fraction of order $p=\P\left[0\in\mathbb{G}_\infty\right]$ of the $\nu$-dimensional section $\Lambda_L^\nu$. The idea is then to work at a larger scale, say $mK$, $m\geq 1$. The interest of using the $mK$ scale is that the process of good $mK$-boxes stochastically dominates a percolation process with parameter as close to one as we wish for $m$ large enough (recall assumption {\bf (SD)}). We follow this strategy at the $mK$ scale and we show in Lemmas \ref{LemmBootstrapingDirSublinearity} and \ref{LemmBoostrapingDirSub2} that it is possible to obtain a good control of the corrector for points in a large fraction of the $mK$-boxes. Finally, we go back to the $K$ scale by finding a $K$-box contained in a suitable $mK$-box from which we can extend the control on the corrector.

In the rest of the section we add the superscripts $K$ and $mK$ to indicate the considered scale.

Let us denote by $e_1,\dots, e_d$ the vectors of the standard basis of $\RR^d$. In order to control the behavior of the corrector at the scale $mK$, for fixed $C, m,\varepsilon$, let us define the mesurable sets:
\begin{align*}
\hspace{-0,2cm}\mathcal{A}_{C,m,\varepsilon}:=\Big\{&\widehat{\xi}\in\widehat{\N}:~\forall i\in \{1,\dots, d\}, \mbox{ for } e=\pm e_i,\forall N\in\NN^*,\\ 
&\mbox{if }j\in\{1,\dots,N\}\mbox{ is s.t. }je\in\mathbb{G}^{mK}_\infty(\widehat{\xi}) \mbox{ then:}\\
&\exists x\in \xi \mbox{ with }\Vor_\xi(x)\cap B^{mK}_{je}\neq\emptyset \mbox{ s.t. } \forall x_0\in \xi \mbox{ s.t. }\Vor_\xi(x_0)\cap B^{mK}_0\neq\emptyset\\
&\big\Vert\chi (\tau_{x_0}\xi, x-x_0)\big\Vert\leq C+\varepsilon N\Big\},
\end{align*}
and 
\begin{align*}
\hspace{-0,2cm}\mathcal{A}_{C,m}:=\Big\{&\widehat{\xi}\in\widehat{\N}:~\forall x,x'\in \xi \mbox{ with }\Vor_\xi(x)\cap B^{mK}_{0}\neq\emptyset,~\Vor_\xi(x')\cap B^{mK}_{0}\neq\emptyset \mbox{ one has:}\\
&\big\Vert\chi (\tau_x\xi, x'-x)\big\Vert\leq C\Big\}.
\end{align*}

For $\nu\in\{1,\dots,d\}$, $n\in\NN^*$ and $\widehat{\xi}\in\widehat{\N}$, let us also define:
\[\Gamma_{n,\nu}^{C, m,\varepsilon}:=\big\{\z\in\Lambda^\nu_n\cap \mathbb{G}_\infty^{mK}:~\tau_{mK\z}\widehat{\xi}\in\mathcal{A}_{C,m}\cap\mathcal{A}_{C,m,\varepsilon}\big\},\]
and
\[\mathbf{\Gamma}_{n}^{C, m,\varepsilon}:=\bigcap_{\nu=1}^d\big\{\z\in\Lambda^d_n:~\z^{(\nu)}\in\Gamma_{n,\nu}^{C,m,\varepsilon}\big\},\]
where $\z^{(\nu)}=(\z_1,\dots,\z_\nu,0,\dots,0)$.

We now prove three intermediary lemmas.
\begin{lemm}\label{LemmDenistyVeryGoodBoxes}
For each $\delta, \varepsilon>0$, there are $C$ and $m$ such that for $\P-a.a.~\widehat{\xi}$, there exists $n_0=n_0(\widehat{\xi}, C, m, \varepsilon,\delta)<\infty$ such that:
\begin{equation}\label{EqLemmDenistyVeryGoodBoxes}
\frac{\#\mathbf{\Gamma}_{n}^{C, m,\varepsilon}}{\#\Lambda^d_n}\geq 1-\delta,\qquad \forall n\geq n_0.
\end{equation}
\end{lemm}
\begin{dem}
Note that, thanks to the union bound, it suffices to show that for any $\nu\in\{1,\dots,d\}$:
\[\frac{\#\Gamma_{n,\nu}^{C, m,\varepsilon}}{\#\Lambda^\nu_n}\geq 1-\frac{\delta}{d}\]
for $n\geq n_0$. Let $\delta':=\delta/d$.

Thanks to assumption {\bf (SD)}, one has for $m$ large enough:
\begin{equation}\label{EqPr1DenistyVeryGoodBoxes}
\P\big[0\in\mathbb{G}^{mK}_{\infty}\big]\geq 1-\frac{\delta'}{2}.
\end{equation}
Given $\delta',\varepsilon$ and $m$ such that the inequality above holds, using Proposition \ref{LemmeSousLineariteDir} at the scale $mK$, we can find $C$ large enough such that:
\begin{equation}\label{EqPr2DenistyVeryGoodBoxes}
\P\big[\mathcal{A}_{C,m}\cap\mathcal{A}_{C,m, \varepsilon}\big\vert 0\in\mathbb{G}^{mK}_\infty\big]\geq 1-\frac{\delta'}{2}.
\end{equation}

Due to the ergodicity assumption {\bf(Er)}, (\ref{EqPr1DenistyVeryGoodBoxes}) and (\ref{EqPr2DenistyVeryGoodBoxes}), one has:
\begin{align*}
\lim_{n\rightarrow\infty}\frac{\#\Gamma_{n,\nu}^{C, m,\varepsilon}}{\#\Lambda^\nu_n}&=\lim_{n\rightarrow\infty}\frac{1}{\#\Lambda^\nu_n}\sum_{\z\in\Lambda_n^\nu}\mathbf{1}_{\tau_{mK\z}\widehat{\xi}\in (\mathcal{A}_{C,m}\cap\mathcal{A}_{C,m, \varepsilon}\cap\{0\in\mathbb{G}_\infty^{mK}\})}\\
&=\P\big[\mathcal{A}_{C,m}\cap\mathcal{A}_{C,m, \varepsilon}\cap\{ 0\in\mathbb{G}^{mK}_\infty\}\big]\\
&=\P\big[\mathcal{A}_{C,m}\cap\mathcal{A}_{C,m, \varepsilon}\big\vert 0\in\mathbb{G}^{mK}_\infty\big]\P\big[0\in\mathbb{G}^{mK}_\infty\big]\\
&\geq \Big(1-\frac{\delta'}{2}\Big)\Big(1-\frac{\delta'}{2}\Big)>1-\delta'.
\end{align*}
This implies the result.
\end{dem}
\renewcommand{\aa}{\mathbf{a}}
\begin{lemm}\label{LemmBootstrapingDirSublinearity}
Given $C,m,\varepsilon >0$ and $\widehat{\xi}\in\widehat{\N}$, if $x\in\xi$ is such that $\Vor_\xi(x)\cap B^{mK}_\aa\neq\emptyset$ for some $\aa\in\mathbf{\Gamma}_n^{C,m,\varepsilon}$, then there exists $x^1\in\xi$ with $\Vor_\xi(x^1)\cap B^{mK}_{\aa^{(1)}}\neq \emptyset$ satisfying:
\begin{equation}\label{EqLemmBootstrapingDirSublinearity}
\big\Vert\chi(\tau_x\xi,x^1-x)\big\Vert\leq (d-1)\big(C+\varepsilon n\big).
\end{equation}
\end{lemm}
\begin{dem}
Since $\aa\in\mathbf{\Gamma}_n^{C,m,\varepsilon}$, $\aa=\aa^{(d)}\in\Gamma_{n,d}^{C,m,\varepsilon}$ and $\aa^{(d-1)}\in\Gamma_{n,d-1}^{C,m,\varepsilon}\subset \mathbb{G}_\infty^{mK}$. In particular, $\tau_{mK\aa}\widehat{\xi}\in\mathcal{A}_{C,m,\varepsilon}$ and we can write $\aa^{(d-1)}=\aa^{(d)}+je_d$, $\vert j\vert\leq n$. This implies that there exists $x^{d-1}\in\xi$ with $\Vor_\xi(x^{d-1})\cap B^{mK}_{\aa^{(d-1)}}\neq\emptyset$ satisfying:
\begin{equation*}
\big\Vert\chi(\tau_x\xi,x^{d-1}-x)\big\Vert\leq C+\varepsilon n.
\end{equation*}
Let us write $x^d:=x$. One can construct in the same way and by induction $x^d,x^{d-1},\dots, x^1$ such that $\Vor_\xi(x^{i-1})\cap B^{mK}_{\aa^{(i-1)}}\neq\emptyset$ and satisfying:
\begin{equation*}
\big\Vert\chi(\tau_{x^i}\xi,x^{i-1}-x^i)\big\Vert\leq C+\varepsilon n, \qquad i=2,\dots,d.
\end{equation*}

The shift-covariance of the corrector leads to:
\[\chi(\tau_x\xi,x^1-x)=\sum_{i=2}^d\left\{\chi(\tau_x\xi,x^{i-1}-x)-\chi(\tau_x\xi,x^i-x)\right\}=\sum_{i=2}^d\chi(\tau_{x^i}\xi,x^{i-1}-x^i).\]

Hence,
\[\big\Vert \chi(\tau_x\xi,x^1-x)\big\Vert\leq\sum_{i=2}^d\big\Vert\chi(\tau_{x^i}\xi,x^{i-1}-x^i)\big\Vert\leq (d-1)\big( C+\varepsilon n\big).\]
\end{dem}
\newcommand{\bb}{\mathbf{b}}
\begin{lemm}\label{LemmBoostrapingDirSub2}
Let  $C,m,\varepsilon >0$, $\widehat{\xi}\in\widehat{\N}$, $x,y\in\xi$ with $\Vor_\xi(x)\cap B^{mK}_\aa\neq\emptyset$ and $\Vor_\xi(y)\cap B^{mK}_\bb\neq\emptyset$ for some $\aa,\bb\in\mathbf{\Gamma}_n^{C,m,\varepsilon}$. Then,
\[\big\Vert\chi (\tau_x\xi,y-x)\big\Vert\leq 2d\big(C+\varepsilon n\big).\]
\end{lemm}
\begin{dem}
Let $x^1,y^1\in\xi$ with $\Vor_\xi(x^1)\cap B^{mK}_{\aa^{(1)}}\neq\emptyset$ and $\Vor_\xi(y^1)\cap B^{mK}_{\bb^{(1)}}\neq\emptyset$ given by Lemma \ref{LemmBootstrapingDirSublinearity}. Thanks to the shift-covariance and the antisymmetry of the corrector, one has:
\begin{align*}
\chi(\tau_x\xi,y-x)
&=\chi(\tau_x\xi,x^1-x)+\chi(\tau_{x^1}\xi,y^1-x^1)+\chi(\tau_{y^1}\xi,y-y^1)\\
&=\chi(\tau_x\xi,x^1-x)+\chi(\tau_{x^1}\xi,y^1-x^1)-\chi(\tau_{y}\xi,y^1-y).
\end{align*}
Together with Lemma \ref{LemmBootstrapingDirSublinearity}, this leads to:
\begin{equation}\label{Eq1ProofBootstraping2}
\big\Vert\chi(\tau_x\xi,y-x)\big\Vert\leq 2(d-1)(C+\varepsilon n) +\big\Vert\chi(\tau_{x^1}\xi,y^1-x^1)\big\Vert.
\end{equation}

Since $\aa^{(1)}\in\Gamma_{n,1}^{C,m,\varepsilon}$ and $\bb^{(1)}\in\mathbb{G}_\infty^{mK}$, there exists a point $\bar{y}\in \xi$ with $\Vor_\xi(\bar{y})\cap B^{mK}_{\bb^{(1)}}\neq\emptyset$ satisfying:
\begin{equation}\label{Eq2ProofBootstraping2}
\big\Vert\chi(\tau_{x^1}\xi,\bar{y}-x^1)\big\Vert\leq C+2\varepsilon n.
\end{equation}

Moreover, $\tau_{mK\bb^{(1)}}\widehat{\xi}\in\mathcal{A}_{C,m}$ which implies that:
\begin{equation}\label{Eq3ProofBootstraping2}
\big\Vert\chi(\tau_{\bar{y}}\xi,y^1-\bar{y})\big\Vert\leq C.
\end{equation}

Bounds (\ref{Eq1ProofBootstraping2})-(\ref{Eq3ProofBootstraping2}) and the shift-covariance of the corrector finally give that:
\[\big\Vert\chi(\tau_x\xi,y-x)\big\Vert\leq 2(d-1)(C+\varepsilon n) +\big\Vert\chi(\tau_{x^1}\xi,\bar{y}-x^1)\big\Vert+\big\Vert\chi(\tau_{\bar{y}}\xi,y^1-\bar{y})\big\Vert\leq 2d( C+\varepsilon n).\]
\end{dem}

\begin{demof}{Proposition \ref{PropSousLineariteMoy}}
We must show that for each $\delta_0, \varepsilon_0>0$, for $\P[\cdot\nolinebreak\vert\nolinebreak 0\nolinebreak\in\nolinebreak\mathbb{G}_\infty]-\nolinebreak a\nolinebreak.\nolinebreak a\nolinebreak.$ $\widehat{\xi}\in\widehat{\N}$, for any $x_0\in\xi$ with $\Vor_\xi(x_0)\cap B^K_0\neq\emptyset$:
\begin{equation*}
\frac{1}{\#\Lambda_L}\sum_{\tiny\begin{array}{c}x\in \xi :\\
\exists \z\in\mathbb{G}^K_\infty\cap\Lambda_L\\
\Vor_\xi(x)\cap B^K_\z\neq\emptyset\end{array}}\mathbf{1}_{\Vert\chi (\tau_{x_0}\xi,x-x_0)\Vert\geq \varepsilon_0L}\leq \delta_0
\end{equation*}
for every $L$ large enough.
Let us define $p:=\P\big[0\in\mathbb{G}_\infty^K\big]$ and  fix $\delta<\min\big(p/2,\delta_0/(2D)\big)$ and $\varepsilon<	\varepsilon_0/(4d)$.
We then choose $C$ and $m$ large enough such that the conclusion of Lemma \ref{LemmDenistyVeryGoodBoxes} holds.
Without  loss of generality we restrict our attention to the case $L=mn$, $n\in\NN$. We will work at both scales $K$ and $mK$.

By ergodicity 
\[\frac{\#\big\{j\in\Lambda^1_L:~je_1\in\mathbb{G}_\infty^{K}\big\}}{2L+1}\xrightarrow[\quad L\longrightarrow \infty\quad]{}p,\]
in particular, for $L$ large enough
\begin{equation}\label{EqProofSubliAver1}\#\big\{j\in\Lambda^1_L:~je_1\in\mathbb{G}_\infty^{K}\big\}\geq\frac{(2L+1)p}{2}.
\end{equation}

On the other hand, denoting by $\pi^{(1)}$ the projection along the first coordinate axis, one has  $\pi^{(1)}\big(\mathbf{\Gamma}^{C,m,\varepsilon}_n\big)\subset \mathbf{\Gamma}^{C,m,\varepsilon}_n$ and 
\[\#\big(\pi^{(1)}\big(\mathbf{\Gamma}^{C,m,\varepsilon}_n\big)\big)\geq\frac{\#\mathbf{\Gamma}^{C,m,\varepsilon}_n}{(2n+1)^{d-1}}\geq(1-\delta)(2n+1)\]
by Lemma \ref{LemmDenistyVeryGoodBoxes}. It follows that
\begin{equation}\label{EqProofSubliAver2}\#\left\{j\in\Lambda^1_L:~\left\lfloor \frac{j}{m}\right\rfloor e_1\in\mathbf{\Gamma}^{C,m,\varepsilon}_n\right\}\geq (1-\delta)(2n+1)m\geq (1-\delta)(2L+1).
\end{equation}

Due to the choice of $\delta$, $(1-\delta)(2L+1)+(2L+1)p/2>(2L+1)=\#\Lambda_L^1$ which implies that 
\[\big\{j\in\Lambda^1_L:~je_1\in\mathbb{G}_\infty^{K}\big\}\cap\left\{j\in\Lambda^1_L:~\left\lfloor \frac{j}{m}\right\rfloor e_1\in\mathbf{\Gamma}^{C,m,\varepsilon}_n\right\}\neq\emptyset.\]

Fix $j$ in the intersection above, thanks to the sublinearity in the direction $e_1$ at scale $K$ (see Proposition \ref{LemmeSousLineariteDir}), there exists $x\in\xi$ with $\Vor_\xi (x)\cap B^K_{je_1}\neq\emptyset$ satisfying:
\[\big\Vert \chi(\tau_{x_0}\xi,x-x_0)\big\Vert\leq C+\varepsilon L, \qquad\forall x_0\in\xi\mbox{ with }\Vor_\xi (x_0)\cap B^K_0\neq\emptyset.\]

Together with Lemma \ref{LemmBoostrapingDirSub2} applied with $\aa=\lfloor \frac{j}{m}\rfloor e_1$ and the shift-covariance of the corrector, this allows us to conclude that for any $y\in\xi$ whose Voronoi cell intersects an $mK$-box with index $\bb$ in $\mathbf{\Gamma}^{C,m,\varepsilon}_n$, for any $x_0\in\xi$ with $\Vor_\xi (x_0)\cap B^K_0\neq\emptyset$:
\begin{align*}
\big\Vert \chi(\tau_{x_0}\xi,y-x_0)\big\Vert&\leq\big\Vert \chi(\tau_{x_0}\xi,x-x_0)\big\Vert+\big\Vert \chi(\tau_{x}\xi,y-x)\big\Vert \\
&\leq C+\varepsilon L+2d\big( C+\varepsilon n\big) \leq (2d+1)\big( C+\varepsilon L\big).
\end{align*}
Thanks to the choice of $\varepsilon$, the last quantity is smaller than $\varepsilon_0L$ when $L$ is large enough. For any $x_0\in \xi$ with $\Vor_\xi(x_0)\cap B^K_0\neq\emptyset$, one has for $L$ large enough:
\begin{align*}
\frac{1}{\#\Lambda_L}\sum_{\tiny\begin{array}{c}x\in \xi :\\
\exists \z\in\mathbb{G}^K_\infty\cap \Lambda_L\\
\Vor_\xi(x)\cap B^K_\z\neq\emptyset\end{array}}\mathbf{1}_{\Vert\chi (\tau_{x_0}\xi,x-x_0)\Vert\geq \varepsilon_0L}
&\leq \frac{1}{\#\Lambda_L}\sum_{\aa\in\Lambda_n\setminus\mathbf{\Gamma}^{C,m,\varepsilon}_n }\sum_{\tiny\begin{array}{c}\z\in\mathbb{G}_\infty^K:\\
B^K_\z\subset B^{mK}_\aa\end{array}}D\\
&\leq \frac{m^dD\#\big(\Lambda_n\setminus\mathbf{\Gamma}^{C,m,\varepsilon}_n\big)}{\#\Lambda_L}\\
&\leq 2D\Bigg(1-\frac{\#\big(\mathbf{\Gamma}^{C,m,\varepsilon}_n\big)}{\#\Lambda_n}\Bigg)\\
&\leq 2D\delta\leq\delta_0
\end{align*}
where we used that:
\[\frac{m^d\#\Lambda_n}{\#\Lambda_L}\xrightarrow[\quad n\rightarrow \infty\quad]{} 1\]
and Lemma \ref{LemmDenistyVeryGoodBoxes}.
\end{demof}
\section{Random walks on $\mathcal{G}_\infty(\widehat{\xi})$}
In order to derive the strong sublinearity of the corrector in $\mathcal{G}_\infty(\widehat{\xi})$ from its sublinearity on average, we need to obtain heat-kernel estimates and bounds on the expected distance between  positions of the walker at time t and 0 (see equations (\ref{EqPropHeatKernelEstimates}) and (\ref{EqPropExpDistBound})). Such estimates cannot be obtained directly for the random walk on the (full) Delaunay triangulation generated by $\xi$ in which the degree is not bounded. Nevertheless, since $\mathcal{G}_\infty(\widehat{\xi})$ has good regularity properties, these bounds will be established in Sections \ref{HeatKernelBound} and \ref{ExpectedDistanceBound} for restricted random walks described below.

For $\widehat{\xi}\in\widehat{\N}$, let us consider the Markov chain $(\widehat{Y}_n)_{n\in\NN}=(\widehat{Y}^{\widehat{\xi}}_n)_{n\in\NN}$ on $\mathcal{G}_\infty(\widehat{\xi})$ induced from the original (discrete time) random walk $(X_n)_{n\in\NN}$ on $\DT (\xi)$. In other words, $(\widehat{Y}_n)_{n\in\NN}$ is the time-homogeneous Markov chain on $\mathcal{G}_\infty(\widehat{\xi})$ with jump probabilities given by:
\begin{equation}\label{EqDefMarcheSautelesTrousDisc}
\widehat{c}^{\widehat{\xi}}_{x,y}:=P^{\widehat{\xi}}\big[\widehat{Y}_{k+1}=y\big\vert\widehat{Y}_k=x\big]=P^{\xi}_x\big[X_{T_1}=y\big], \qquad x,y\in\mathcal{G}_\infty ({\widehat{\xi}}),
\end{equation}
where $T_1:=\inf\{j\geq 1:~X_j\in\mathcal{G}_\infty (\widehat{\xi})\}$. Note that the \emph{holes} (\emph{i.e.} the ($\DT(\xi)$-)connected components of $\xi\setminus\mathcal{G}_\infty (\widehat{\xi})$) are a.s. finite. Hence, $T_1$ is a.s. finite and $(\widehat{Y}_n)_{n\in\NN}$ is well defined. Moreover, for any $z\in\mathcal{G}_\infty (\widehat{\xi})$, applying the optional stopping theorem to the martingale $(X_n-z-\chi (\tau_z\xi, X_n-z))_{n\in\NN}$ starting from $z$, it appears that, $(\widehat{Y}_n-z-\chi (\tau_z\xi, \widehat{Y}_n-z))_{n\in\NN}$ is a martingale.

We also consider a continuous-time version of the random walk defined above $(\widehat{Y}_t)_{t\geq 0}:=(\widehat{Y}_{N(t)})_{t\geq 0}$ where $N(t)$ is the intensity 1 Poisson process on the half-line $\RR^+$. It has infinitesimal generator:
\begin{equation}\label{EqDefMarcheSautelesTrousGene}
\widehat{\mathcal{L}}^{\widehat{\xi}}f(x):=\sum_{y\in\mathcal{G}_{\infty}(\widehat{\xi})}\widehat{c}^{\widehat{\xi}}_{x,y}\big( f(y)-f(x)\big), \qquad x,y\in\mathcal{G}_\infty ({\widehat{\xi}}).
\end{equation}

It is not difficult to see that $(\widehat{Y}_t-z-\chi (\tau_z\xi, \widehat{Y}_t-z))_{t\geq 0}$ is also a martingale.

We denote by $P^{\widehat{\xi}}_x$ the (quenched) law of this walk starting from $x$. Note that this walk has speed `at most 1'. Observe also that the measure $\deg_{\DT (\xi)}$ is reversible w.r.t. both $(\widehat{Y}_n)_{n\in\NN}$ and $(\widehat{Y}_t)_{t\geq 0}$. Actually, standard computations show that the \emph{detailed balance condition}:
\[\deg_{\DT(\xi)}(x)\widehat{c}^{\widehat{\xi}}_{x,y}=\deg_{\DT(\xi)}(y)\widehat{c}^{\widehat{\xi}}_{y,x}, \qquad x,y\in\mathcal{G}_\infty(\widehat{\xi}),\]
is satisfied.
\section{Heat-kernel estimates for $(\widehat{Y}^{\widehat{\xi}}_t)_{t>0}$}\label{HeatKernelBound}
The aim of this section is to prove the following heat-kernel bound.
\begin{prop}\label{PropHeatKernelEstimates}
For a.a. $\widehat{\xi}\in\widehat{\N}$:
\begin{equation}\label{EqPropHeatKernelEstimates}
\sup_{n\geq 1}\max_{x\in\mathcal{G}_\infty(\widehat{\xi})\cap [-n,n]^d}\sup_{t\geq n}t^{\frac{d}{2}}P^{\widehat{\xi}}_x\left[\widehat{Y}_t=x\right]<\infty
\end{equation}
where $(\widehat{Y}_t)_{t\geq 0}$ is the continuous-time random walk on $\mathcal{G}_\infty(\widehat{\xi})$  with generator (\ref{EqDefMarcheSautelesTrousGene}).
\end{prop}

The proof of this bound relies on isoperimetric inequalities and the technics developed in \cite{MoPe}; it is completed in Subsection \ref{SousSectProofofHKE}. Precise definitions are given in Subsection \ref{SousSectionPrecDef}. Isoperimetric inequalities for random walks confined in large boxes are established in Subsection \ref{SousSectIso}. Additional technical results are isolated from the proof of Proposition \ref{PropHeatKernelEstimates} and given in Subsection \ref{SousSectTechnicalLemmas}.
\subsection{Precise definitions}\label{SousSectionPrecDef}
We will state isoperimetric inequalities for random walks confined in large boxes. We need to introduce additional notations and random walks confined in boxes of side $L$.

Recall the definitions of $\mathbb{G}_{(L)}$  and $\mathcal{G}_{(L)}(\widehat{\xi})$ from the introduction and denote by $\overline{\mathbb{G}}_{(L)}$ the complementary in $\ZZ^d$ of the unique unbounded ($l^1-$)connected component of $\ZZ^d\setminus \mathbb{G}_{(L)}$ and 
by $\overline{\mathcal{G}}_{(L)}(\widehat{\xi})$ the set of points of $\xi$ whose Voronoi cell intersects a $K$-box with index in  $\overline{\mathbb{G}}_{(L)}$. In other words, $\overline{\mathbb{G}}_{(L)}$ is the union of $\mathbb{G}_{(L)}$ with the (discrete) holes contained in $[-L,L]^d$ and  $\overline{\mathcal{G}}_{(L)}(\widehat{\xi})$ is the union of $\mathcal{G}_{(L)}(\widehat{\xi})$ with the holes (for the $\DT(\xi)$ structure) contained in $\left[ -(L+\frac{1}{2})K,(L+\frac{1}{2})K\right]^d$. Then, we write $\left(\overline{X}^{(L)}_n\right)_{n\in\NN}$ for the random walk in the restriction of $\DT(\xi)$ to $\overline{\mathcal{G}}_{(L)}(\widehat{\xi})$. Denoting by $T^*_n$ the time of $n^\mathrm{th}$ visit to $\mathcal{G}_{(L)}(\widehat{\xi})$ for this walk, we consider the induced discrete-time random walk $\left(\widehat{Y}^{(L)}_n\right)_{n\in\NN}:=\left(\overline{X}^{(L)}_{T^*_n}\right)_{n\in\NN}$ on $\mathcal{G}_{(L)}(\widehat{\xi})$. We also consider its continuous-time counterpart $\left(\widehat{Y}^{(L)}_t\right)_{t\geq 0}:=\left(\widehat{Y}^{(L)}_{N(t)}\right)_{t\geq 0}$ where $N(\cdot)$ is an intensity 1 Poisson process on the half-line $\RR^+$. The random walk $\left(\widehat{Y}^{(L)}_t\right)_{t\geq 0}$ can be thought of as the continuous-time random walk on $\mathcal{G}_{(L)}(\widehat{\xi})$ with speed at most 1 and which `jumps holes' of $\mathcal{G}_{(L)}(\widehat{\xi})$. Let us denote by $\widehat{\deg}_{L,\widehat{\xi}}(\cdot)$ the degree in the restriction of $\DT(\xi)$ to $\overline{\mathcal{G}}_{(L)}(\widehat{\xi})$ and write $\widehat{\deg}_{L,\widehat{\xi}}(A):=\sum_{x\in A}\widehat{\deg}_{L,\widehat{\xi}}(x)$, $A\subset\overline{\mathcal{G}}_{(L)}(\widehat{\xi})$. Classical computations show that $\widehat{\deg}_{L,\widehat{\xi}}(\cdot)$ is reversible for both $\left(\overline{X}^{(L)}_n\right)_{n\in\NN}$ and $\left(\widehat{Y}^{(L)}_t\right)_{t\geq 0}$. The \emph{conductance} of the set $A\subset \mathcal{G}_{(L)}(\widehat{\xi})$ w.r.t. $\left(\widehat{Y}^{(L)}_t\right)_{t\geq 0}$ is given by:
\begin{equation}\label{DefConduct1}
\widehat{I}^{(L)}_A=\widehat{I}^{\widehat{\xi}, (L)}_A:=\frac{\sum_{x\in A}\sum_{y\in\mathcal{G}_{(L)}(\widehat{\xi})\setminus A}\widehat{\deg}_{L,\widehat{\xi}}(x)\overline{P}^{\widehat{\xi},(L)}_x\left[\overline{X}^{(L)}_{T^*_1}=y\right]}{\widehat{\deg}_{L,\widehat{\xi}}(A)},
\end{equation}
where $\overline{P}^{\widehat{\xi},(L)}_x$ stands for the law of $\left(\overline{X}^{(L)}_n\right)_{n\in\NN}$. The associated \emph{isoperimetric profile} is:
\begin{equation}\label{EqDefIsoProf1}
\widehat{\varphi}_{(L)}(u):=\inf\left\{\widehat{I}^{(L)}_A:~\widehat{\deg}_{L,\widehat{\xi}}(A)\leq u\widehat{\deg}_{L,\widehat{\xi}}\left(\mathcal{G}_{(L)}(\widehat{\xi})\right)\right\}, \qquad u\in\left]0,\frac{1}{2}\right].
\end{equation}

The advantage of $\left(\widehat{Y}^{(L)}_t\right)_{t\geq 0}$ is that this walk coincides with $\left(\widehat{Y}_t\right)_{t\geq 0}$ as long as they do not leave (the interior of) $\overline{\mathcal{G}}_{(L)}(\widehat{\xi})$. Nevertheless, we are not able to obtain directly a bound on the isoperimetric profile $\widehat{\varphi}_{(L)}$. In a similar way as in \cite{BBHK}, we compare it with the isoperimetric profile of the constant speed random walk on $\mathcal{G}_{(L)}(\widehat{\xi})$, that is the walk $\left(\widetilde{Y}^{(L)}_t\right)_{t\geq 0}$ with generator:
\begin{equation}\label{EqDefGeneMarcheRestreinte}
\widetilde{\mathcal{L}}_{\widehat{\xi}}^{(L)}(x,y):=\left\{\begin{array}{ll}
\dfrac{\mathbf{1}_{x\sim y\text{ in }\DT(\xi)}}{\widetilde{\deg}_{L,\widehat{\xi}}(x)}&\mbox{if }x\neq y\\
-1 &\mbox{if }x= y
\end{array}\right. ,
\end{equation} 
where $\widetilde{\deg}_{L,\widehat{\xi}}(\cdot )$ denotes the degree in the restriction of $\DT(\xi)$ to $\mathcal{G}_{(L)}(\widehat{\xi})$.
The measure $\widetilde{\deg}_{L,\widehat{\xi}}(\cdot)$ is clearly reversible w.r.t. $\left(\widetilde{Y}^{(L)}_t\right)_{t\geq 0}$. The associated \emph{conductance} of the set $A\subset \mathcal{G}_{(L)}(\widehat{\xi})$ is given by:
\begin{equation}\label{DefConduct2}
\widetilde{I}^{(L)}_A=\widetilde{I}^{\widehat{\xi}, (L)}_A:=\frac{\sum_{x\in A}\sum_{y\in\mathcal{G}_{(L)}(\widehat{\xi})\setminus A}\mathbf{1}_{x\sim y\text{ in }\DT(\xi)}}{\widetilde{\deg}_{L,\widehat{\xi}}(A)},
\end{equation}
and the corresponding \emph{isoperimetric profile} is:
\begin{equation}\label{EqDefIsoProf2}
\widetilde{\varphi}_{(L)}(u):=\inf\left\{\widetilde{I}^{(L)}_A:~\widetilde{\deg}_{L,\widehat{\xi}}(A)\leq u\widetilde{\deg}_{L,\widehat{\xi}}\left(\mathcal{G}_{(L)}(\widehat{\xi})\right)\right\}, \qquad u\in\left]0,\frac{1}{2}\right].
\end{equation}
\subsection{Isoperimetric inequality}\label{SousSectIso}
The goal of this section is to obtain a lower  bound on the isoperimetric profile $\widehat{\varphi}_{(L)}$; it is stated in Corollary \ref{IsoHat}.
\subsubsection{Comparison between $\widetilde{\varphi}_{(L)}$ and $\widehat{\varphi}_{(L)}$}
First, note that for $x\in\mathcal{G}_{(L)}(\widehat{\xi})$:
\begin{equation}\label{EqCompIsoProf1}
D\geq\deg_{\DT (\xi)}(x)\geq\widehat{\deg}_{L,\widehat{\xi}}(x)\geq\widetilde{\deg}_{L,\widehat{\xi}}(x)\geq 1\geq\frac{\widehat{\deg}_{L,\widehat{\xi}}(x)}{D}\geq\frac{\widetilde{\deg}_{L,\widehat{\xi}}(x)}{D},
\end{equation}
and that for $x,y\in\mathcal{G}_{(L)}(\widehat{\xi})$
\begin{equation}
\overline{P}^{\widehat{\xi},(L)}_x\left[\overline{X}^{(L)}_{T^*_1}=y\right]\geq \overline{P}^{\widehat{\xi},(L)}_x\left[\overline{X}^{(L)}_{1}=y\right]\mathbf{1}_{x\sim y\text{ in }\DT(\xi)}=\frac{\mathbf{1}_{x\sim y\text{ in }\DT(\xi)}}{\widehat{\deg}_{L,\widehat{\xi}}(x)}. 
\end{equation}

Hence, for $A\subset\mathcal{G}_{(L)}(\widehat{\xi})$, we have:
\[\sum_{x\in A}\sum_{y\in\mathcal{G}_{(L)}(\widehat{\xi})\setminus A} \widehat{\deg}_{L,\widehat{\xi}}(x) \overline{P}^{\widehat{\xi},(L)}_x\left[\overline{X}^{(L)}_{T^*_1}=y\right]\geq \sum_{x\in A}\sum_{y\in\mathcal{G}_{(L)}(\widehat{\xi})\setminus A}\mathbf{1}_{x\sim y\text{ in }\DT(\xi)}.\]
With $\widehat{\deg}_{L,\widehat{\xi}}(A)\leq D \widetilde{\deg}_{L,\widehat{\xi}}(A)$, this implies that:
\begin{equation}\label{EqCompConductances}
\widehat{I}^{(L)}_A\geq \frac{\widetilde{I}^{(L)}_A}{D}.
\end{equation}

Using (\ref{EqCompIsoProf1}) and (\ref{EqCompConductances}), one deduces:
%
%
\begin{lemm}\label{LemmCompIsoProfile}
For $u\in\left]0,\frac{1}{2D}\right]$:
\begin{equation}\label{EqLemmCompIsoProfile}
\widehat{\varphi}_{(L)}(u)\geq\frac{\widetilde{\varphi}_{(L)}(Du)}{D}.
\end{equation}
\end{lemm}
\subsubsection{Lower bounds for $\widetilde{\varphi}_{(L)}$}
Our aim is to show the following bound on the isoperimetric profile $\widetilde{\varphi}_{(L)}$ associated with $\left(\widetilde{Y}^{(L)}_t\right)_{t\geq 0}$.
\begin{prop}\label{IsoTilde} There exists $c=c (d,K,\alpha)>0$ such that $\mathcal{P}-a.s.$ for $L$ large enough:
\[\widetilde{\varphi}_{(L)}(u)\geq c\min\Bigg\{\dfrac{1}{u^{1/d}L},\dfrac{1}{\log(L)^{\frac{d}{d-1}}}\Bigg\}, \quad u\in \Big]0,\frac{1}{2}\Big].\]
\end{prop}

As in \cite{CFIso}, we use as much as possible an isoperimetric inequality for the percolation cluster $\mathbb{G}_{(L)}$.

\begin{prop}[see \cite{CFIso}, eq. (2.5)]
\label{IsoBP}
There exists $\kappa>0$ such that almost surely for $L$ large enough, for $\mathbb{A}\subset\mathbb{G}_{(L)}$ with $0<\#(\mathbb{A})\leq\frac{1}{2}\#(\mathbb{G}_{(L)})$:
\[\dfrac{\#(\partial \mathbb{A})}{\#(\mathbb{A})}\geq\kappa \min\Bigg\{\dfrac{1}{\#(\mathbb{A})^{\frac{1}{d}}},\dfrac{1}{\log(L)^{\frac{d}{d-1}}}\Bigg\},\]
where $\partial \mathbb{A}=\{x\in\mathbb{G}_{(L)}\setminus \mathbb{A}:\,x\sim y\mbox{ for some }y\in \mathbb{A}\}$ is the (vertex external) boundary of the set $\mathbb{A}$ in $\mathbb{G}_{(L)}$.
\end{prop}

This result can be proved by adapting the arguments given in \cite[Appendix]{BBHK} to the context of supercritical site percolation (see also the proof of \cite[Lemma 2.6]{BenjaminiMossel} for $p \gg p_c^{\text{site}}	(\ZZ^d)$). 

We adapt the proof of \cite[Theorem 1.1]{CFIso} to the present setting. It is worth noting that the arguments of \cite{CFIso} can be used to derive isoperimetric bounds for the Delaunay triangulation confined in cubic boxes at least when the underlying point process is a PPP. This does not lead to sharp enough heat-kernel bounds for the random walk on the full Delaunay triangulation due to the unboundedness of the degree.

\bigskip
\begin{demof}{Proposition \ref{IsoTilde}}
For $A\subset\mathcal{G}_{(L)}(\widehat{\xi})$, we define:
\[\mathbb{L}(A):=\big\{\z\in\mathbb{G}_{(L)}:\,\exists x\in A\mbox{ s.t. }\operatorname{Vor}_\xi(x)\cap B_\z\neq\emptyset\big\}.\] 
Let us observe that, thanks to the definition of the good boxes, for any $A\subset\mathcal{G}_{(L)}(\widehat{\xi})$:
\begin{equation}\label{EqCompVol}
\frac{\# \mathbb{L}(A)}{2^d}\leq \widetilde{\deg}_{L,\widehat{\xi}}(A)\leq \# A\max_{x\in A}\deg_{\DT(\xi)}(x)\leq D^2\# \mathbb{L}(A). 
\end{equation}
In the first inequality, we used that, for $x\in\mathcal{G}_{(L)}(\widehat{\xi})$, $\operatorname{Vor}_{\xi}(x)$ does not intersect more than $2^d$ good boxes since the diameter of the cell is less than $K$.

From now on we assume that $\widetilde{\deg}_{L,\widehat{\xi}}(A)\leq\frac{1}{2}\widetilde{\deg}_{L,\widehat{\xi}}(\mathcal{G}_{(L)}(\widehat{\xi}))$. We are going to discuss separately the cases when $\#(\mathbb{L}(A))$ is large or small with respect to $\#(\mathbb{G}_{(L)})$. Roughly, if $\#(\mathbb{L}(A))$ is large then $\#(\mathbb{L}(\mathcal{G}_{(L)}(\widehat{\xi})\setminus A))$ is not too small and $\widetilde{I}^{(L)}_A$ is easily bounded from below by some constant. When $\#(\mathbb{L}(A))$ is small, a bound is obtained using the isoperimetric inequality for $\mathbb{G}_{(L)}$ given in Proposition \ref{IsoBP}.

\paragraph{{\it The case $\#(\mathbb{L}(A))>\big(1-\frac{1}{2^{d+2}D^2}\big)\#(\mathbb{G}_{(L)})$}.}
Using the general bound (\ref{EqCompVol}) and inequality $\widetilde{\deg}_{L,\widehat{\xi}}(A)\leq\frac{1}{2}\widetilde{\deg}_{L,\widehat{\xi}}(\mathcal{G}_{(L)}(\widehat{\xi}))$, one obtains:
\[\#\mathbb{L}(\mathcal{G}_{(L)}(\widehat{\xi})\setminus A)\geq \frac{\widetilde{\deg}_{L,\widehat{\xi}}(\mathcal{G}_{(L)}(\widehat{\xi})\setminus A)}{D^2}\geq \frac{\widetilde{\deg}_{L,\widehat{\xi}}(\mathcal{G}_{(L)}(\widehat{\xi}))}{2D^2}\geq \frac{\#\mathbb{G}_{(L)}}{2^{d+1}D^2}.
\]

It follows that $\mathbb{L}(A)$ and $\mathbb{L}(\mathcal{G}_{(L)}(\widehat{\xi})\setminus A)$ have a large intersection in this case:
\begin{align}\label{EqLargeIntersection}
\#\big(\mathbb{L}(A)\cap\mathbb{L}(\mathcal{G}_{(L)}(\widehat{\xi})\setminus A)\big)&\geq \# \mathbb{L}(A)-\#\big(\mathbb{G}_{(L)}\setminus\mathbb{L}(\mathcal{G}_{(L)}(\widehat{\xi})\setminus A)\big) \nonumber\\
&= \#\mathbb{L}(A)-\#(\mathbb{G}_{(L)})+\#\big(\mathbb{L}(\mathcal{G}_{(L)}(\widehat{\xi})\setminus A)\big) \nonumber\\
&\geq\Big(1-\frac{1}{2^{d+2}D^2}\Big)\#(\mathbb{G}_{(L)})-\#(\mathbb{G}_{(L)})+\frac{1}{2^{d+1}D^2}\#(\mathbb{G}_{(L)}) \nonumber\\
&=\frac{1}{2^{d+2}D^2}\#(\mathbb{G}_{(L)}).
\end{align}
This allows us to bound from below the numerator in $\widetilde{I}^{(L)}_A$ as follows. If $\z\in \mathbb{L}(A)\cap\mathbb{L}(\mathcal{G}_{(L)}(\widehat{\xi})\setminus A)$, one can choose $x\in A$ and $y\in\mathcal{G}_{(L)}(\widehat{\xi})\setminus A$  whose respective Voronoi cells intersect $B_{\z}$. Hence, connecting the nuclei of the Voronoi cells which intersect the line segment $[x,y]$, it is easy  to find an edge between a point of $A$ and a point of $\mathcal{G}_{(L)}(\widehat{\xi})\setminus A$ which is included in $\overline{B}_\z=\bigcup_{\z':\Vert \z'-\z\Vert_\infty\leq 1}B_{\z'}$ thanks to the definition of good boxes. Since a specific edge is associated to at most $3^d$ boxes by this procedure, it follows using (\ref{EqLargeIntersection}) that:
\begin{equation*}\label{EqBoundNumCondL}
\sum_{x\in A}\sum_{y\in \mathcal{G}_{(L)}(\widehat{\xi})\setminus A}\mathbf{1}_{x\sim y\text{ in }\DT(\xi)}\geq \frac{\#\big(\mathbb{L}(A)\cap\mathbb{L}(\mathcal{G}_{(L)}(\widehat{\xi})\setminus A)\big)}{3^d}\geq\frac{\#(\mathbb{G}_{(L)})}{4\cdot 6^dD^2}.
\end{equation*} 

Since $\widetilde{\deg}_{L,\widehat{\xi}}(A)\leq D^2\# \big(\mathbb{L}(A)\big)\leq D^2\#(\mathbb{G}_{(L)})$, we obtain that:
\begin{equation*}
\widetilde{I}^{(L)}_A\geq \frac{1}{4\cdot 6^dD^4}.
\end{equation*} 
\paragraph{{\it The case $\#(\mathbb{L}(A))\leq\big(1-\frac{1}{2^{d+2}D^2}\big)\#(\mathbb{G}_{(L)})$}.}
Let us show that, in this case, the numerator in $\widetilde{I}^{(L)}_A$ can be bounded from below in terms of $\#\partial \mathbb{L}(A)$ or $\#\partial (\mathbb{G}_{(L)}\setminus\mathbb{L}(A))$. If $B_\z$ and $B_{\z'}$ are two neighboring good boxes such that $\z\in \mathbb{L}(A)$ and $\z'\in\mathbb{G}_{(L)}\setminus\mathbb{L}(A)$, there exists an edge between a point of $A$ and a point of $\mathcal{G}_{(L)}(\widehat{\xi})\setminus A$ contained in $\overline{B}_\z\cup \overline{B}_{\z'}$. To see this, let us fix a point $x\in A$ whose Voronoi cell intersects $B_{\z}$ and a point $y\in\mathcal{G}_{(L)}(\widehat{\xi})\setminus A$ whose Voronoi cell intersects $B_{\z'}$. It then suffices to connect the consecutive nuclei of the Voronoi cells which intersect the line segment $[x,y]$ to find an edge between a point of $A$ and a point of $\mathcal{G}_{(L)}(\widehat{\xi})\setminus A$. This edge is contained in $\overline{B}_\z\cup \overline{B}_{\z'}$ thanks to the definition of good boxes. It follows that there exists $\delta=\delta(d)$ such that:
\begin{equation*}
\sum_{x\in A}\sum_{y\in \mathcal{G}_{(L)}(\widehat{\xi})\setminus A}\mathbf{1}_{x\sim y\text{ in }\DT(\xi)}\geq \delta\max\big\{\#\partial \mathbb{L}(A),\#\partial (\mathbb{G}_{(L)}\setminus\mathbb{L}(A))\big\}.
\end{equation*} 

Since $\widetilde{\deg}_{L,\widehat{\xi}}(A)\leq D^2\# \mathbb{L}(A)\leq D^2(2^{d+2}D^2-1)\# \big(\mathbb{G}_{(L)}\setminus\mathbb{L}(A)\big)$, we deduce that:
\begin{equation*}
\widetilde{I}^{(L)}_A\geq \frac{\delta}{(2^{d+2}D^2-1)D^2}\frac{\#\partial\mathbb{A}}{\#\mathbb{A}},
\end{equation*}
for $\mathbb{A}=\mathbb{L}(A)\mbox{ and }\mathbb{A}=\mathbb{G}_{(L)}\setminus\mathbb{L}(A)$.

Choosing
\[\mathbb{A}:=\left\{
\begin{array}{ll}
\mathbb{L}(A) &\mbox{if }\#\mathbb{L}(A)\leq\frac{1}{2}\#\mathbb{G}_{(L)}\\
\mathbb{G}_{(L)}\setminus\mathbb{L}(A)&\mbox{otherwise}
\end{array}
\right. , \]
Proposition \ref{IsoBP} then implies that almost surely for $L$ large:
\begin{align*}
\widetilde{I}^{(L)}_A&\geq\frac{\kappa\delta}{(2^{d+2}D^2-1)D^2}\min\Bigg\{\dfrac{1}{\#(\mathbb{A})^{\frac{1}{d}}},\dfrac{1}{\log(L)^{\frac{d}{d-1}}}\Bigg\}\\
&\geq\frac{\kappa\delta}{(2^{d+2}D^2-1)D^2}\min\Bigg\{\dfrac{1}{\#(\mathbb{L}(A))^{\frac{1}{d}}},\dfrac{1}{\log(L)^{\frac{d}{d-1}}}\Bigg\}.
\end{align*}

Using that $\# \mathbb{L}(A)\leq 2^d\widetilde{\deg}_{L,\widehat{\xi}}(A)$, we obtain that:
\begin{equation*}
\widetilde{I}^{(L)}_A\geq\frac{\kappa\delta}{2(2^{d+2}D^2-1)D^2}\min\Bigg\{\dfrac{1}{\widetilde{\deg}_{L,\widehat{\xi}}(A)^{\frac{1}{d}}},\dfrac{1}{\log(L)^{\frac{d}{d-1}}}\Bigg\}.
\end{equation*}

Since $\widetilde{\deg}_{L,\widehat{\xi}} (\mathcal{G}_{(L)}(\widehat{\xi}))\leq D^2\#\mathbb{G}_{(L)}\leq D^2(2L+1)^d$,
the conclusion of Proposition \ref{IsoTilde} follows.
\end{demof}

Combining Lemma \ref{LemmCompIsoProfile} and Proposition \ref{IsoTilde}, we obtain:
\begin{coro}\label{IsoHat}
There exists $c=c (d,K,\alpha)>0$ such that $\mathcal{P}-a.s.$ for $L$ large enough:
\[\widehat{\varphi}_{(L)}(u)\geq c\min\Bigg\{\dfrac{1}{u^{1/d}L},\dfrac{1}{\log(L)^{\frac{d}{d-1}}}\Bigg\}, \quad u\in \Big]0,\frac{1}{2D}\Big].\]
\end{coro} 
\subsection{Other technical results}\label{SousSectTechnicalLemmas}
\subsubsection{Volume growth for $\widehat{\deg}_{L,\widehat{\xi}}\left(\mathcal{G}_{(L)}(\widehat{\xi})\right)$}
We briefly check that there exist constants $c$ and $C$ such that a.s. for $L$ large enough:
\begin{equation}\label{EqCroissVolHat}
cL^d\leq \widehat{\deg}_{L,\widehat{\xi}}\left(\mathcal{G}_{(L)}(\widehat{\xi})\right)\leq C L^d.
\end{equation}
The upper bound is very simple since it suffices to write:
\begin{align*}
\widehat{\deg}_{L,\widehat{\xi}}\left(\mathcal{G}_{(L)}(\widehat{\xi})\right)&\leq D \#\left(\mathcal{G}_{(L)}(\widehat{\xi})\right)\leq D\max_{\z\in\mathbb{G}_{(L)}}\#\left\{x\in\xi:~\Vor_\xi (x)\cap B_\z\neq\emptyset\right\}\#\left(\mathbb{G}_{(L)}\right)\\
&\leq D^2 \#\left([-L,L]^d\cap\ZZ^d\right)\leq CL^d.
\end{align*}

Since any box with index in $\mathbb{G}_{(L)}$ contains at least a point of $\mathcal{G}_{(L)}(\widehat{\xi})$ which has degree at least 1:
\[\widehat{\deg}_{L,\widehat{\xi}}\left(\mathcal{G}_{(L)}(\widehat{\xi})\right)\geq \#\left(\mathcal{G}_{(L)}(\widehat{\xi})\right)\geq\#\left(\mathbb{G}_{(L)}\right). \]
The lower bound follows by using that a.s. for $L$ large enough:
\[\#\left(\mathbb{G}_{(L)}\right)\geq c L^d,\]
which is a consequence of the ergodic theorem.
\subsubsection{Size of the holes and connectivity of $\mathcal{G}_\infty(\widehat{\xi})$ in large boxes}
In order to compare $(\widehat{Y}_t)_{t\geq 0}$ with $(\widehat{Y}^{(L)}_t)_{t\geq 0}$, we need to control the size of the holes ({\it i.e.} $\DT(\xi)$-connected components of $\xi\setminus\mathcal{G}_\infty(\widehat{\xi})$) and to establish connectivity properties of $\mathcal{G}_\infty(\widehat{\xi})$ in large boxes. More precisely, for $C, \gamma>1$ and $t>0$, let define the events: 
\begin{align*}
\mathcal{A}_t=\mathcal{A}_{t,\gamma, C}:=\Bigg\{&\mbox{any hole contained in }\left[-K\left(\lfloor t^\gamma\rfloor+\frac{1}{2}\right),K\left(\lfloor t^\gamma\rfloor+\frac{1}{2}\right)\right]^d\\
&\mbox{ has diameter smaller than }C\log t\Bigg\},
\end{align*}
 and
\begin{align*}
\mathcal{B}_t=\mathcal{B}_{t,\gamma, C}&:=\left\{\mathcal{G}_\infty(\widehat{\xi})\cap Q_t
=\mathcal{G}_{(\lfloor t^\gamma\rfloor)}(\widehat{\xi})\cap Q_t\right\}
\end{align*}
where:
\begin{align*}
Q_t&:=\Bigg[-t\left(\frac{3C}{2}\log t+1\right),t\left(\frac{3C}{2}\log t+1\right)\Bigg]^d.
\end{align*}
We prove that:
\begin{lemm}\label{LemmSizeHolesAndConnectivity}
Assuming that $K$ is large enough, for each $\gamma>1$, there exists $C<\infty$ such that almost surely for $t$ large enough $\mathcal{A}_t$ and $\mathcal{B}_t$ are realized. 
\end{lemm}
\begin{dem} First, observe that any $\DT (\xi)$-connected component $A$ of $\xi\setminus\mathcal{G}_\infty (\widehat{\xi})$ is contained in the union of $K$-boxes with indices in some discrete hole $\mathbb{A}$ ({\it i.e.} a connected component of $\ZZ^d\setminus\mathbb{G}_\infty$). Hence, in order to show that $\mathcal{A}_t$ holds almost surely, it suffices to verify that a.s. any discrete hole contained in $\left[-\lfloor t^\gamma\rfloor,\lfloor t^\gamma\rfloor\right]^d\cap \ZZ^d$ has diameter at most $C\log t$, for $C$ suitably chosen. Denote by $\mathbb{A}_\z$ the (possibly empty) hole at $\z\in\ZZ^d$ for an independent percolation process of parameter $p$. Recall that assumption {\bf (SD)} ensures that the process of `good boxes' dominates such a process with $p$ as close to 1 as we wish whenever $K$ is fixed large enough. Assuming that $p$ is large enough, a standard Peierls argument shows that there exists $c_6>0$ such that:
\[\P\left[\diam \mathbb{A}_0\geq n\right]\leq e^{-c_6n}.\]
Thus, 
\[\P\left[\exists \z\in [-n,n]^d\cap\ZZ^d\mbox{ s.t. } \diam \mathbb{A}_\z\geq C \log n\right]\leq (2n+1)^d e^{-c_6C\log n}.\]
Our first claim then follows by the Borel-Cantelli lemma if $C$ is well chosen.  

As above, in order to show that $\mathcal{B}_t$ is a.s. realized for $t$ large enough, we only need to check the corresponding claim for the percolation process, that is: almost surely for $t$ large enough, $\mathbb{G}_\infty\cap Q_t=\mathbb{G}_{(\lfloor t^\gamma\rfloor)}\cap Q_t$.

Let us justify that, almost surely for $t$ large enough, $\mathbb{G}_{(\lfloor t^\gamma\rfloor)}$ coincides with the largest connected component  of $\mathbb{G}_{\infty}\cap \left[-\lfloor t^\gamma\rfloor,\lfloor t^\gamma\rfloor\right]^d$. Thanks to \cite[Lemma B.1]{CFIso} and the Borel-Cantelli lemma, we know that there exists $L_0$ a.s. finite such that, for $L\geq L_0$, the maximal open cluster $\mathbb{G}_{(L)}$ in $[-L,L]^d\cap \ZZ^d$  is the only open cluster in this box with diameter larger than $L/10$ and crosses this box in every coordinate direction (see also \cite[Remark 7]{CFIso}). In particular, $\mathbb{G}_{(L)}$ has diameter $2L\geq (L+1)/10$ and is thus included in $\mathbb{G}_{(L+1)}$. So, for $L'\geq L\geq L_0$, $\mathbb{G}_{(L)}$ is contained in $\mathbb{G}_{(L')}$. Hence, it is included in an open cluster with infinite diameter which is $\mathbb{G}_\infty$.

It remains to verify that any two vertices $\z$ and $\z'$ of $\mathbb{G}_{\infty}\cap Q_t$ are connected by an open path whithin $\left[-\lfloor t^\gamma\rfloor,\lfloor t^\gamma\rfloor\right]^d$. Call $\mathcal{B}'_t$ this event and write $\d_{\mathbb{G}_\infty}(\cdot, \cdot )$ for the graph distance in $\mathbb{G}_\infty$. Assume that $\mathcal{B}'_t$ fails for some large $t$ and fix $\z, \z'\in\mathbb{G}_{\infty}\cap Q_t$ which are not connected in  $\left[-\lfloor t^\gamma\rfloor,\lfloor t^\gamma\rfloor\right]^d$. Considering a shortest path from $\z$ to $\z'$ one can find $\z''$ in $\mathbb{G}_{\infty}\cap Q_{3t}\setminus Q_{2t}$ such that $\d_{\mathbb{G}_\infty}(\z', \z'' )\geq \lfloor t^\gamma\rfloor-t(3C\log (t)/2+1)$. Hence, for any $\kappa$, for $t$ large enough:
\begin{equation*}\label{EqProofSizeHoleCon1}
\P\left[ (\mathcal{B}'_t)^c\right]\leq \sum_{\z'\in Q_t}\sum_{\z''\in Q_{3t}\setminus Q_{2t}}\P\left[\z',\z''\in\mathbb{G}_\infty,~\d_{\mathbb{G}_\infty}(\z', \z'' )\geq \kappa \Vert \z''-\z'\Vert\right].
\end{equation*}

From \cite{AP}, we know that $\P\left[\z',\z''\in\mathbb{G}_\infty,~\d_{\mathbb{G}_\infty}(\z', \z'' )\geq \kappa \Vert \z''-\z'\Vert\right]$ decreases exponentially with $\Vert \z''-\z'\Vert$ which implies that $\P\left[ (\mathcal{B}'_t)^c\right]$ decreases exponentially with $t$. One finally concludes thanks to the Borel-Cantelli lemma.
\end{dem}
\subsection{Proof of Proposition \ref{PropHeatKernelEstimates}}\label{SousSectProofofHKE}
For $t\geq 0$, let us denote by $\widehat{N}(t)$ the number of jumps of $(\widehat{Y}_s)_{s\geq 0}$ up to time $t$ and by $\mathcal{C}_t$ the event: `$\widehat{N}(t)\leq 3t/2$'. Since  $(\widehat{Y}_s)_{s\geq 0}$ has speed at most 1, $\widehat{N}(\cdot)$ is dominated by a Poisson process of intensity 1 on $\RR^+$ and $P^{\widehat{\xi}}_x\left[\mathcal{C}_t^c\right]\leq c_7\exp(-c_8t)$ for some $c_7,c_8>0$. This implies that, almost surely for $t$ large enough, $\mathcal{C}_t$ is realized and we only need to obtain the heat-kernel bound on this event.

Recall the definitions of $\mathcal{A}_t$ and $\mathcal{B}_t$ from the previous subsection. On $\mathcal{A}_t\cap\mathcal{B}_t\cap\mathcal{C}_t$, starting from a point of $\mathcal{G}_\infty(\widehat{\xi})\cap [-t,t]^d$, $(\widehat{Y}_s)_{s\geq 0}$ does at most $3/2t$ jumps of length at most $C\log t$ up to time $t$. In particular, it has visited only points of $\mathcal{G}_{(\lfloor t^\gamma\rfloor)}(\widehat{\xi})\cap Q_t$ and does not depend on $\xi\setminus\overline{\mathcal{G}}_{(\lfloor t^\gamma\rfloor)}(\widehat{\xi})$ up to time $t$. Hence, we can find a coupling of $(\widehat{Y}_s)_{s\geq 0}$ and $(\widehat{Y}^{(\lfloor t^\gamma\rfloor)}_s)_{s\geq 0}$ such that these two coincide up to time $t$. 
Thus, for $t$ large enough, we can write:
\begin{align*}
P^{\widehat{\xi}}_{x}\left[\widehat{Y}_t=y\right]&=P^{\widehat{\xi}}_{x}\left[\left\{\widehat{Y}_t=y\right\}\cap \mathcal{A}_t\cap\mathcal{B}_t\cap\mathcal{C}_t\right]\\
&=P^{\widehat{\xi}, (\lfloor t^\gamma\rfloor )}_{x}\left[\left\{\widehat{Y}^{(\lfloor t^\gamma\rfloor )}_t=y\right\}\cap \mathcal{A}_t\cap\mathcal{B}_t\cap\mathcal{C}_t\right]\\
&= P^{\widehat{\xi}, (\lfloor t^\gamma\rfloor )}_{x}\left[\widehat{Y}^{(\lfloor t^\gamma\rfloor )}_t=y\right].
\end{align*}
It then remains to bound $P^{\widehat{\xi}, (\lfloor t^\gamma\rfloor )}_{x}\left[\widehat{Y}^{(\lfloor t^\gamma\rfloor )}_t=y\right]$. To this end, we will rely on the isoperimetric inequality stated in Corollary \ref{IsoHat} and apply the strategy developed by Morris and Peres in \cite{MoPe}.
\begin{theo}[see {\cite[Theorem 13]{MoPe}}]\label{ThMoPe}
Let $(X_t)_{t\geq 0}$ be an irreducible continuous-time Markov chain on a finite state space $\mathcal{X}$ with reversible probability measure $\pi$ and isoperimetric profile $\varphi$.  
 
 For all $\varepsilon>0$ and all $x,y\in\mathcal{X}$, if 
\begin{equation}\label{EqCondMoPe}
t\geq \int_{4\min (\pi(x),\pi(y))}^{4/\varepsilon}\frac{8 \d u}{u\varphi(u)^2}
\end{equation}
then
\[P_{x}\left[X_t=y\right]\leq \pi(y )\left(1+\varepsilon\right).\]
\end{theo}
Here the reversible probability measure is given by:
\[\widehat{\pi}_{\lfloor t^\gamma\rfloor}(\cdot ):=\widehat{\deg}_{\lfloor t^\gamma\rfloor, \widehat{\xi} }(\cdot)/\widehat{\deg}_{\lfloor t^\gamma\rfloor, \widehat{\xi} }(\mathcal{G}_{(\lfloor t^\gamma\rfloor)}(\widehat{\xi})).\]
Choosing $\varepsilon$ of the form $\varepsilon=c_9t^{d(\gamma-\frac{1}{2})}$ ($c_9$ will be chosen large enough), the conclusion of Theorem \ref{ThMoPe} reads:
\[P^{\widehat{\xi}, (\lfloor t^\gamma\rfloor )}_{x}\left[\widehat{Y}^{(\lfloor t^\gamma\rfloor )}_t=y\right]\leq \frac{\widehat{\deg}_{\lfloor t^\gamma\rfloor, \widehat{\xi} }( y)}{\widehat{\deg}_{\lfloor t^\gamma\rfloor, \widehat{\xi} }(\mathcal{G}_{(\lfloor t^\gamma\rfloor)}(\widehat{\xi}))}\left(1+c_9t^{d(\gamma-\frac{1}{2})}\right).\]
Using Lemma \ref{LemmBoundedDegree} and (\ref{EqCroissVolHat}), one deduces that $P^{\widehat{\xi}, (\lfloor t^\gamma\rfloor )}_{x}\left[\widehat{Y}^{(\lfloor t^\gamma\rfloor )}_t=y\right]\leq c_{10}t^{-\frac{d}{2}}$ for $t$ large enough. It remains to check the validity of (\ref{EqCondMoPe}) with $\varepsilon=c_9t^{d(\gamma-\frac{1}{2})}$ for $t$ large enough when $c_9$ is well chosen. 

Assuming that $t$ satisfies $4t^{d(\frac{1}{2}-\gamma)}\leq c_9/(2D)$, one obtains with Corollary \ref{IsoHat} and inequality (\ref{EqCompVol}) that:
\begin{align*}
& \int_{4\min (\widehat{\pi}_{\lfloor t^\gamma\rfloor }(x),\widehat{\pi}_{\lfloor t^\gamma\rfloor }(y))}^{4c_9^{-1}t^{d(\frac{1}{2}-\gamma)}}\frac{8 \d u}{u\widehat{\varphi}_{(\lfloor t^\gamma\rfloor)}(u)^2}\\
&\qquad\qquad\qquad\leq 8 \int_{c_{11}/t^{\gamma d}}^{4c_9^{-1}t^{d(\frac{1}{2}-\gamma)}}\frac{\d u}{u\widehat{\varphi}_{(\lfloor t^\gamma\rfloor)}(u)^2}\\
&\qquad\qquad\qquad\leq c_{12}\left((\log t)^{\frac{2d}{d-1}}\int_{c_{11}/t^{\gamma d}}^{(\log \lfloor t^\gamma\rfloor)^\frac{d^2}{d-1}/\lfloor t^\gamma\rfloor^d}\frac{\d u}{u}+t^{2\gamma} \int_{(\log \lfloor t^\gamma\rfloor)^\frac{d^2}{d-1}/\lfloor t^\gamma\rfloor^d}^{4c_9^{-1}t^{d(\frac{1}{2}-\gamma)}}u^{\frac{2}{d}-1}\d u
\right)\\
 &\qquad\qquad\qquad\leq c_{12}\left(\frac{d}{2}\left(\frac{4}{c_9}\right)^\frac{2}{d}t+c_{13}(\log t)^{\frac{2d}{d-1}}\left(\log\log t +c_{14}\right)\right).
\end{align*}

If $c_9$ has been fixed large enough, the last expression is smaller than $t$ for every $t$ large enough.

To summarize, we have just proved that for a.a. $\widehat{\xi}$, there exist $c_{15}=c_{15}(\widehat{\xi})$ and $T=T(\widehat{\xi})$ such that for any $t\geq T$, for any $x\in\mathcal{G}_\infty(\widehat{\xi})\cap [-t,t]^d$ and any $y\in\mathcal{G}_\infty(\widehat{\xi})$:
\[P^{\widehat{\xi}}_{x}\left[\widehat{Y}_t=y\right]\leq \frac{c_{15}}{t^\frac{d}{2}}.\]  

This implies the required result.\begin{flushright}
$\square$
\end{flushright}
\section{Expected distance bound for $(\widehat{Y}^{\widehat{\xi}}_t)_{t>0}$}\label{ExpectedDistanceBound}
It is known that bounds on the expected distance between the position of the walk at time $t$ and its starting point can be derived from the heat-kernel estimate (\ref{EqPropHeatKernelEstimates}) as soon as the volume grows regularly (see for example \cite{Bass, Barlow, BP}).
In this section, we use this strategy to prove the following proposition.
\begin{prop}\label{PropExpDistBound}
For a.a. $\widehat{\xi}\in\widehat{\N}$:
\begin{equation}\label{EqPropExpDistBound}
\sup_{n\geq 1}\max_{x\in\mathcal{G}_\infty(\widehat{\xi})\cap [-n,n]^d}\sup_{t\geq n}\frac{E^{\widehat{\xi}}_x\left[\left\Vert\widehat{Y}_t-x\right\Vert\right]}{\sqrt{t}}<\infty.
\end{equation}
\end{prop} 
\begin{dem}
Proposition \ref{PropHeatKernelEstimates} shows that the assumption of \cite[Proposition 6.2]{BP} is satisfied in the present setting. Hence, there exist constants $c_{16}$ and $c_{17}$ such that for a.a. $\widehat{\xi}$, for every $x\in\mathcal{G}_{\infty}(\widehat{\xi})$, for $t$ large enough:
\begin{align}\label{EqProofExpDisBound1}\frac{E^{\widehat{\xi}}_x\left[\widehat{\d}\left(\widehat{Y}_t,x\right)\right]}{\sqrt{t}}&\leq c_{16}+c_{17}\sup_{0<s\leq t^{-\frac{1}{2}}}\left\{s^d\sum_{y\in \mathcal{G}_{\infty}(\widehat{\xi})}e^{-s\widehat{\d}(x,y)}\right\},
\end{align}
where $\widehat{\d}\left(x,y\right)=\widehat{\d}_{\mathcal{G}_{\infty}(\widehat{\xi})}\left(x,y\right)$ denotes the `natural' distance between $x$ and $y$ for  $(\widehat{Y}_t)_{t>0}$ (\emph{i.e.} the minimal number of jumps that the random walk needs to do in order to go from $x$ to $y$).

At this point, we need to compare $\widehat{\d}$ with the Euclidean distance and to check that the r.h.s. of (\ref{EqProofExpDisBound1}) is uniformly bounded. Let us denote by $\widehat{\d}'=\widehat{\d}'_{\mathbb{G}_{\infty}}$ the chemical distance in $\mathbb{G}_{\infty}$ in which we add an edge between every two points on the boundary of a (shared) discrete hole. For $x\in\mathcal{G}_\infty (\widehat{\xi})$, we choose $\z (x)\in\mathbb{G}_\infty$ such that $\Vor_\xi(x)$ intersects $B_{\z(x)}$ to be the minimal one in the lexicographic order. It is not difficult to see that the definitions of $\mathbb{G}_\infty$ and $\mathcal{G}_\infty (\widehat{\xi})$ imply that there are constants $c_{18}$ and $c_{19}$ such that:
\begin{equation}\label{EqProofExpDisBound2}
c_{18}\widehat{\d}'(\z(x),\z(y))\leq\widehat{\d}(x,y)\leq c_{19}(\widehat{\d}'(\z(x),\z(y))+1), \qquad \forall x,y\in\mathcal{G}_\infty (\widehat{\xi}).
\end{equation}
By the same arguments as in the proof of \cite[Lemma 3.1]{BP}, one obtains that:
\[\P\left[\z,\z'\in\mathbb{G}_\infty,~\widehat{\d}'(\z,\z')\leq c_{20} \Vert\z-\z'\Vert\right]\leq e^{-c_{21} \Vert\z-\z'\Vert},\]
for suitable constants $c_{20}$ and $c_{21}$. 

Using the estimate above and the Borel-Cantelli lemma, one deduces that there exists a constant $C$ such that, for $n\geq N=N(\widehat{\xi})$, for all $\z\in\mathbb{G}_\infty\cap [-n/K,n/K]^d$, for all $\z'\in\mathbb{G}_\infty$ with $\Vert \z'-\z\Vert\geq C \log n$:
\begin{equation}
\label{EqProofExpDistBound3}\widehat{\d}'(\z,\z')\geq c_{20} \Vert \z'-\z\Vert.
\end{equation}  
Then, (\ref{EqProofExpDisBound1}), (\ref{EqProofExpDisBound2}) and (\ref{EqProofExpDistBound3}) imply that for every $x\in\mathcal{G}_{\infty}(\widehat{\xi})\cap [-n,n]^d$, for $t\geq n \geq N$:
\begin{equation}\label{EqProofExpDistBound4}
\frac{E^{\widehat{\xi}}_x\left[\left\Vert\widehat{Y}_t-x\right\Vert\right]}{\sqrt{t}}\leq c_{22}+c_{23}\sup_{0<s\leq t^{-\frac{1}{2}}}\left\{s^d\sum_{y\in \mathcal{G}_{\infty}(\widehat{\xi})}e^{-s\widehat{\d}(x,y)}\right\}.
\end{equation}

Finally, observe that once $1/s\geq t^\frac{1}{2}\geq n^\frac{1}{2}\gg C\log n$ for any $x\in \mathcal{G}_\infty(\widehat{\xi})\cap [-n,n]^d$:
\begin{align*}
\sum_{y\in \mathcal{G}_{\infty}(\widehat{\xi})}e^{-s\widehat{\d}(x,y)}&\leq \sum_{\z'\in \mathbb{G_\infty}}\sum_{\tiny\begin{array}{c}y\in \xi\\ \Vor_\xi(y)\cap B_{\z'}\neq\emptyset\end{array}}e^{-s\widehat{\d}(x,y)}\\
&\leq D\sum_{\z'\in \mathbb{G_\infty}}e^{-c_{18}s\widehat{\d}'(\z(x),\z')}\\
&\leq c_{24}\left(s^{-d}+\sum_{\tiny\begin{array}{c}\z'\in \mathbb{G}_\infty:\\ \Vert\z'-\z(x)\Vert\geq 1/s\end{array}}e^{-c_{25}s\Vert \z'-\z(x)\Vert}\right)\\
&\leq c_{26}s^{-d}.
 \end{align*}

The conclusion then follows thanks to (\ref{EqProofExpDistBound4}).
\end{dem}
\section{Almost sure sublinearity in $\mathcal{G}_\infty(\widehat{\xi})$}\label{ASRestSubli}
The aim of this section is to prove the `strong' sublinearity of the corrector in $\mathcal{G}_\infty (\widehat{\xi})$.
\begin{prop}\label{PropStrogSubliGoodPoints}
For $\P[\cdot\vert 0\in \mathbb{G}_\infty]-$a.a. $\widehat{\xi}$:
\begin{equation}\label{EqPropStrogSubliGoodPoints}
\lim_{n\rightarrow\infty}\frac{\mathfrak{R}_n}{n}=0,
\end{equation}
where
\[\mathfrak{R}_n=\mathfrak{R}_n(\widehat{\xi}):=\max_{\tiny\begin{array}{c}x_0\in \xi :\\\Vor_\xi(x_0)\cap B^K_0\neq\emptyset\end{array}}\max_{x\in\mathcal{G}_\infty(\widehat{\xi})\cap [-n,n]^d}\left\Vert \chi (\tau_{x_0}\xi,x-x_0)\right\Vert.\]
\end{prop}

Following an idea attributed to Yuval Peres in \cite{BP} and \cite{CFP}, it suffices to show the recursive bound:
\begin{lemm}\label{LemmRecStrogSubliGoodPoints}
For $\P[\cdot\vert 0\in \mathbb{G}_\infty]-$a.e. $\widehat{\xi}$, for each $\varepsilon,\delta>0$, there exists $n_0=n_0(\widehat{\xi},\varepsilon,\delta)<\infty$ such that:
\begin{equation}\label{EqLemmRecStrogSubliGoodPoints}
\mathfrak{R}_n\leq \varepsilon n+\delta \mathfrak{R}_{3n},\qquad \forall n\geq n_0.
\end{equation}
\end{lemm}

For the reader's convenience, we recall how Proposition \ref{PropStrogSubliGoodPoints} can be deduced from the previous lemma and Proposition \ref{LemmCroissPol} (see also \cite[Proof of Theorem 2.4, p. 1337]{BP}). Assume that the conclusion of Proposition \ref{PropStrogSubliGoodPoints} is false and choose $0<c_{27}<\limsup_{n\rightarrow\infty}\mathfrak{R}_n/n$, $\varepsilon:=c_{27}/2$ and $\delta:=1/3^{\beta+1}$ with $\beta$ such that:
\begin{equation}\label{EqProofStrongSubliGoodPoints1}
\lim_{n\rightarrow\infty}\frac{\mathfrak{R}_n}{n^\beta}=0.
\end{equation}

For infinitely many $n$'s, $\mathfrak{R}_n\geq c_{27}n$ which implies by Lemma \ref{LemmRecStrogSubliGoodPoints} that:
\[\mathfrak{R}_{3n}\geq\frac{\mathfrak{R}_n-\varepsilon n}{\delta}\geq \frac{(c_{27}-\varepsilon)n}{\delta}\geq 3^\beta c_{27} n\]
when $n\geq n_0$. One then obtains by induction that $\mathfrak{R}_{3^kn}\geq c_{27}3^{\beta k}n$ which contradicts (\ref{EqProofStrongSubliGoodPoints1}). We now turn our attention to the proof of Lemma \ref{LemmRecStrogSubliGoodPoints}.

\begin{demof}{Lemma \ref{LemmRecStrogSubliGoodPoints}}
We adapt the arguments given in \cite[\S 5]{BP}.
Let us define:
\[C_1=C_1(\widehat{\xi}):=\sup_{n\geq 1}\max_{x\in\mathcal{G}_\infty(\widehat{\xi})\cap [-n,n]^d}\sup_{t\geq n}t^{\frac{d}{2}}P^{\widehat{\xi}}_x\left[\widehat{Y}_t=x\right]\]
and 
\[C_2=C_2(\widehat{\xi}):=\sup_{n\geq 1}\max_{x\in\mathcal{G}_\infty(\widehat{\xi})\cap [-n,n]^d}\sup_{t\geq n}\frac{E^{\widehat{\xi}}_x\left[\left\Vert\widehat{Y}_t-x\right\Vert\right]}{\sqrt{t}}.\]
Recall that these two quantities are a.s. finite thanks to Propositions \ref{PropHeatKernelEstimates} and \ref{PropExpDistBound}.

For large $n$, we choose $y_0=y_0(n)$ with $\Vor_\xi (y_0)\cap B_0\neq \emptyset$ and $y=y(n)\in [-n,n]^d$ such that $\mathfrak{R}_n=\Vert \chi(\tau_{y_0}\xi,y-y_0)\Vert$ and we define the stopping time:
\[S_n:=\inf\left\{t\geq 0:~\left\Vert \widehat{Y}_t-y\right\Vert\geq 2n\right\}.\]

For $n$ is large enough, holes have sizes of logarithmic order (see Lemma \ref{LemmSizeHolesAndConnectivity}) and thus $\Vert \widehat{Y}_{t\wedge S_n}-y\Vert\leq 3n$ for all $t$. Due to the harmonicity of $\varphi$, the optional stopping theorem gives:
\[E^{\widehat{\xi}}_y\left[\widehat{Y}_{t\wedge S_n}-y-\chi \left(\tau_y\xi,\widehat{Y}_{t\wedge S_n}-y\right)\right]=E^{\widehat{\xi}}_y\left[\varphi\left( \tau_y\xi,\widehat{Y}_{t\wedge S_n}-y\right)\right]=0.\]

By the shift-covariance of the corrector, one has:
\[\chi \left(\tau_{y_0}\xi,y-y_0\right)=\chi \left(\tau_{y_0}\xi,\widehat{Y}_{t\wedge S_n}-y_0\right)-\chi \left(\tau_y\xi,\widehat{Y}_{t\wedge S_n}-y\right),\]
thus
\[\chi \left(\tau_{y_0}\xi,y-y_0\right)=E^{\widehat{\xi}}_y\left[\chi \left(\tau_{y_0}\xi,\widehat{Y}_{t\wedge S_n}-y_0\right)-\widehat{Y}_{t\wedge S_n}+y\right].\]
It follows that:
\begin{equation} \label{EqProofLemmRec1}
\mathfrak{R}_n=\left\Vert\chi \left(\tau_{y_0}\xi,y-y_0\right)\right\Vert\leq E^{\widehat{\xi}}_y\left[\left\Vert\chi \left(\tau_{y_0}\xi,\widehat{Y}_{t\wedge S_n}-y_0\right)-\widehat{Y}_{t\wedge S_n}+y\right\Vert\right].
\end{equation}

Let us fix $\varepsilon>0$ and define:
\[\mathcal{O}_n:=\left\{x\in\mathcal{G}_\infty(\widehat{\xi})\cap[-3n,3n]^d:\exists x_0\mbox{ s.t.}\Vor_\xi(x_0)\cap B_0\neq\emptyset\mbox{ and}\left\Vert\chi \left(\tau_{x_0}\xi,x-x_0\right)\right\Vert\geq\frac{\varepsilon}{2}n\right\}.\]
Note that, by Proposition \ref{PropSousLineariteMoy}, $\#\mathcal{O}_n=o(n^d)$. Restricting our attention to $t=t(n)\geq 4n$ (whose value will be specified at the end), we will decompose the expectation above as:
\[E^{\widehat{\xi}}_y\left[\left\Vert\chi \left(\tau_{y_0}\xi,\widehat{Y}_{t\wedge S_n}-y_0\right)-\widehat{Y}_{t\wedge S_n}+y\right\Vert\right]=E_1+E_2,\]
with
\[E_1=E^{\widehat{\xi}}_y\left[\left\Vert\chi \left(\tau_{y_0}\xi,\widehat{Y}_{t\wedge S_n}-y_0\right)-\widehat{Y}_{t\wedge S_n}+y\right\Vert\mathbf{1}_{S_n< t}\right],\]
and
\[E_2=E^{\widehat{\xi}}_y\left[\left\Vert\chi \left(\tau_{y_0}\xi,\widehat{Y}_{t\wedge S_n}-y_0\right)-\widehat{Y}_{t\wedge S_n}+y\right\Vert\mathbf{1}_{S_n\geq t}\right].\]

We first deal with the term $E_1$. Since $t\geq 4n$, Markov inequality shows that:
\[P^{\widehat{\xi}}_y\left[\left\Vert\widehat{Y}_{2t}-y\right\Vert\geq\frac{3}{2}n\right]\leq \frac{2E^{\widehat{\xi}}_y\left[\left\Vert\widehat{Y}_{2t}-y\right\Vert\right]}{3n}\leq\frac{2\sqrt{2t}C_2}{3n}.\]
Observe that $\left\{\left\Vert\widehat{Y}_{2t}-y\right\Vert\leq 3n/2,S_n<t\right\}\subset \left\{\left\Vert\widehat{Y}_{2t}-\widehat{Y}_{S_n}\right\Vert\geq n/2,S_n<t\right\}$. On $\left\{S_n<t\right\}$,  since $s:=2t-S_n\in [t,2t]$, one has:
\[P^{\widehat{\xi}}_z\left[\left\Vert\widehat{Y}_{s}-z\right\Vert\geq\frac{1}{2}n\right]\leq \frac{2\sqrt{2t}C_2}{n},\]
with $z=\widehat{Y}_{S_n}$, this implies by the strong Markov property that:
\[P^{\widehat{\xi}}_y\left[\left\Vert\widehat{Y}_{2t}-y\right\Vert\leq\frac{3}{2}n,S_n<t\right]\leq \frac{2\sqrt{2t}C_2}{n}.\]
Recall that $\left\Vert\widehat{Y}_{t\wedge S_n}-z\right\Vert\leq 3n$ for $n$ large enough. It follows that: 
\begin{align}\label{EqProofLemmRec2}
&E^{\widehat{\xi}}_y\left[\left\Vert\chi \left(\tau_{y_0}\xi,\widehat{Y}_{t\wedge S_n}-y_0\right)-\widehat{Y}_{t\wedge S_n}+y\right\Vert\mathbf{1}_{S_n<t}\right]\nonumber\\
&\qquad\qquad\qquad\leq \left(\mathfrak{R}_{3n}+3n\right)P^{\widehat{\xi}}_y\left[S_n<t\right]\nonumber\\
&\qquad\qquad\qquad\leq \left(\mathfrak{R}_{3n}+3n\right)\left(P^{\widehat{\xi}}_y\left[\left\Vert\widehat{Y}_{2t}-y\right\Vert\geq\frac{3}{2}n\right]+P^{\widehat{\xi}}_y\left[\left\Vert\widehat{Y}_{2t}-y\right\Vert\leq\frac{3}{2}n,S_n<t\right]\right)\nonumber\\
&\qquad\qquad\qquad\leq \frac{8\sqrt{2t}C_2}{3n}\left(\mathfrak{R}_{3n}+3n\right).
\end{align}

Thanks to the definitions of $C_2$ and $\mathcal{O}_n$ and to the fact that $t\geq n$, one has:
\begin{align}\label{EqProofLemmRec3}
&E^{\widehat{\xi}}_y\left[\left\Vert\chi \left(\tau_{y_0}\xi,\widehat{Y}_{t\wedge S_n}-y_0\right)-\widehat{Y}_{t\wedge S_n}+y\right\Vert\mathbf{1}_{S_n\geq t}\right]\nonumber\\
&\quad\leq E^{\widehat{\xi}}_y\left[\left\Vert\widehat{Y}_{t}-y\right\Vert\mathbf{1}_{S_n\geq t}\right]+E^{\widehat{\xi}}_y\left[\left\Vert\chi \left(\tau_{y_0}\xi,\widehat{Y}_{t}-y_0\right)\right\Vert\mathbf{1}_{S_n\geq t}\right]\nonumber\\
&\quad\leq C_2\sqrt{t}+E^{\widehat{\xi}}_y\left[\left\Vert\chi \left(\tau_{y_0}\xi,\widehat{Y}_{t}-y_0\right)\right\Vert\mathbf{1}_{S_n\geq t,\widehat{Y}_{t}\not\in\mathcal{O}_n}\right]\nonumber\\
&\qquad\qquad\qquad\qquad+E^{\widehat{\xi}}_y\left[\left\Vert\chi \left(\tau_{y_0}\xi,\widehat{Y}_{t}-y_0\right)\right\Vert\mathbf{1}_{S_n\geq t,\widehat{Y}_{t}\in\mathcal{O}_n}\right]\nonumber\\
&\quad\leq C_2\sqrt{t}+\frac{\varepsilon}{2}n+\mathfrak{R}_{3n}P^{\widehat{\xi}}_y\left[\widehat{Y}_{t}\in\mathcal{O}_n\right]=C_2\sqrt{t}+\frac{\varepsilon}{2}n+\mathfrak{R}_{3n}\sum_{z\in\mathcal{O}_n}P^{\widehat{\xi}}_y\left[\widehat{Y}_{t}=z\right].
\end{align}

But, using that $\deg_{\DT(\xi)}$ is reversible w.r.t $\left(\widehat{Y}_s\right)_{s\geq 0}$ and bounded by $D$ on $\mathcal{G}_\infty(\widehat{\xi})$, we obtain by the Markov property and the Cauchy-Schwarz inequality that:
\begin{align}\label{EqProofLemmRec4}
P^{\widehat{\xi}}_y\left[\widehat{Y}_{t}=z\right]^2&=\left(\sum_{x\in\mathcal{G}_\infty(\widehat{\xi})}P^{\widehat{\xi}}_y\left[\widehat{Y}_{\frac{t}{2}}=x\right]P^{\widehat{\xi}}_x\left[\widehat{Y}_{\frac{t}{2}}=z\right]\right)^2\nonumber\\
&\leq\left(\sum_{x\in\mathcal{G}_\infty(\widehat{\xi})}P^{\widehat{\xi}}_y\left[\widehat{Y}_{\frac{t}{2}}=x\right]^2\right)\left(\sum_{x\in\mathcal{G}_\infty(\widehat{\xi})}P^{\widehat{\xi}}_x\left[\widehat{Y}_{\frac{t}{2}}=z\right]^2\right)\nonumber\\
&\leq\left(\sum_{x\in\mathcal{G}_\infty(\widehat{\xi})}P^{\widehat{\xi}}_y\left[\widehat{Y}_{\frac{t}{2}}=x\right]\frac{\deg_{\DT(\xi)}(x)}{\deg_{\DT(\xi)}(y)}P^{\widehat{\xi}}_x\left[\widehat{Y}_{\frac{t}{2}}=y\right]\right)\nonumber\\
&\qquad\qquad\times\left(\sum_{x\in\mathcal{G}_\infty(\widehat{\xi})}\frac{\deg_{\DT(\xi)}(z)}{\deg_{\DT(\xi)}(x)}P^{\widehat{\xi}}_z\left[\widehat{Y}_{\frac{t}{2}}=x\right]P^{\widehat{\xi}}_x\left[\widehat{Y}_{\frac{t}{2}}=z\right]\right)\nonumber\\
&\leq D^2P^{\widehat{\xi}}_y\left[\widehat{Y}_{t}=y\right]P^{\widehat{\xi}}_z\left[\widehat{Y}_{t}=z\right]\leq \frac{D^2C_1^2}{t^d}. 
\end{align}

Combining bounds (\ref{EqProofLemmRec1})-(\ref{EqProofLemmRec4}), we get:
\begin{align*}
\mathfrak{R}_n&\leq\frac{8\sqrt{2t}C_2}{3n}\left(\mathfrak{R}_{3n}+3n\right)+C_2\sqrt{t}+\frac{\varepsilon}{2}n+\frac{DC_1\#\mathcal{O}_n}{t^\frac{d}{2}}\mathfrak{R}_{3n}\\
&=\left(\frac{8\sqrt{2t}C_2}{3n}+\frac{DC_1\#\mathcal{O}_n}{t^\frac{d}{2}}\right)\mathfrak{R}_{3n}+\left(\frac{8\sqrt{2t}C_2}{n}+\frac{\varepsilon}{2}\right)n+C_2\sqrt{t}.
\end{align*}

The conclusion then follows by choosing $t=c_{28}n^2$ for some $c_{28}=c_{28}(\widehat{\xi},\varepsilon,\delta)$ small enough and using that $\#\mathcal{O}_n=o(n^d)$.
\end{demof}
\section{Proof of main results}\label{SubliProof}
\subsection{Proof of Theorem \ref{TheoSousLine}}
We first show that Proposition \ref{PropStrogSubliGoodPoints} and the control of the diameter of the holes imply that for $\P\left[\cdot\vert0\in\mathbb{G}_\infty\right]$-a.a. $\widehat{\xi}$:
\begin{equation}\label{EqProofThPrinc1}
\lim_{n\to \infty}\frac{1}{n}\max_{\tiny\begin{array}{c}x_0\in \xi :\\\Vor_\xi(x_0)\cap B^K_0\neq\emptyset\end{array}}\max_{x\in\xi\cap [-n,n]^d}\left\Vert \chi (\tau_{x_0}\xi,x-x_0)\right\Vert=0.
\end{equation} 

Recall that, almost surely for $n$ large enough, holes intersecting $[-n,n]^d$ have diameters smaller than $C\log n$ (see Lemma \ref{LemmSizeHolesAndConnectivity}).
Let $\mathcal{H}\subset\xi$ be a hole intersecting $[-n,n]^d$ and denote by $\partial_{\text{ext}}\mathcal{H}\subset\mathcal{G}_\infty (\widehat{\xi})$ its external boundary, that is the set of the points of $\xi\setminus\mathcal{H}$ which are neighbors of a point of $\mathcal{H}$ in $\DT (\xi)$. We can assume that $\partial_{\text{ext}}\mathcal{H}$ is contained in $[-2n,2n]^d$. Let us define $S:=\inf\{k\geq 0:X_k\not\in\mathcal{H}\}$. Thanks to the harmonicity of $\varphi$ and the optional stopping theorem, for $x\in\mathcal{H}$, one has:
\[E^\xi_x\left[X_S-x-\chi(\tau_x\xi,X_S-x)\right]=0.\] 
But the shift-covariance of the corrector implies that for any $x_0\in\xi$:
\[\chi(\tau_x\xi,X_S-x)=\chi(\tau_{x_0}\xi,X_S-x_0)-\chi(\tau_{x_0}\xi,x-x_0),\]
and thus 
\[\chi(\tau_{x_0}\xi,x-x_0)=E^\xi_x\left[X_S-x-\chi(\tau_{x_0}\xi,X_S-x_0)\right].\]
It follows that:
\begin{align*}
\max_{\tiny\begin{array}{c}x_0\in \xi :\\\Vor_\xi(x_0)\cap B^K_0\neq\emptyset\end{array}}&\max_{x\in\xi\cap [-n,n]^d}\left\Vert \chi (\tau_{x_0}\xi,x-x_0)\right\Vert\\
&\leq \max_{\tiny\begin{array}{c}x_0\in \xi :\\\Vor_\xi(x_0)\cap B^K_0\neq\emptyset\end{array}}\max_{x\in\mathcal{G}_\infty(\widehat{\xi})\cap [-2n,2n]^d}\left\Vert \chi (\tau_{x_0}\xi,x-x_0)\right\Vert+C\log n.
\end{align*}
Together with Proposition \ref{PropStrogSubliGoodPoints}, this implies that (\ref{EqProofThPrinc1}) holds for $\P\left[\cdot\vert0\in\mathbb{G}_\infty\right]$-a.a. $\widehat{\xi}$.

We now prove that (\ref{EqProofThPrinc1}) actually holds for $\P$-a.a. $\xi$. Note that this is the step where we eliminate the coupling. Observe that:
\begin{align*}
\max_{\tiny\begin{array}{c}x_0\in \tau_{Ke_1}\xi :\\\Vor_{\tau_{Ke_1}\xi}(x_0)\cap B^K_0\neq\emptyset\end{array}}&\max_{x\in\tau_{Ke_1}\xi\cap [-n,n]^d}\left\Vert \chi (\tau_{x_0}\tau_{Ke_1}\xi,x-x_0)\right\Vert\\
&=\max_{\tiny\begin{array}{c}x_0\in \xi :\\\Vor_{\xi}(x_0)\cap B^K_{-e_1}\neq\emptyset\end{array}}\max_{x\in\xi\cap  [-n-K,n-K]\times[-n,n]^{d-1}}\left\Vert \chi (\tau_{x_0}\xi,x-x_0)\right\Vert\\
&\leq\max_{\tiny\begin{array}{c}x_1\in \xi :\\\Vor_{\xi}(x_1)\cap B^K_{0}\neq\emptyset\end{array}}\max_{x\in\xi\cap  [-2n,2n]^{d}}\left\Vert \chi (\tau_{x_1}\xi,x-x_1)\right\Vert\\
&\qquad\qquad+\max_{\tiny\begin{array}{c}x_0\in \xi :\\\Vor_{\xi}(x_0)\cap B^K_{-e_1}\neq\emptyset\end{array}}\max_{\tiny\begin{array}{c}x_1\in \xi :\\\Vor_{\xi}(x_1)\cap B^K_{0}\neq\emptyset\end{array}}\left\Vert \chi (\tau_{x_1}\xi,x_0-x_1)\right\Vert.\\
\end{align*}
Hence, the event:
\[\mathcal{A}:=\left\{\lim_{n\to \infty}\frac{1}{n}\max_{\tiny\begin{array}{c}x_0\in \xi :\\\Vor_\xi(x_0)\cap B^K_0\neq\emptyset\end{array}}\max_{x\in\xi\cap [-n,n]^d}\left\Vert \chi (\tau_{x_0}\xi,x-x_0)\right\Vert=0\right\}\]
is shift invariant w.r.t. $\tau=\tau^{K,e_1}$. Thanks to the ergodicity assumption {\bf (Er)}, it is a 0-1 event and we already know that $\P[\mathcal{A}]\geq\P[0\in\mathbb{G}_\infty]>0$. Thus, (\ref{EqProofThPrinc1}) holds $\P$-a.s..

In particular, (\ref{EqProofThPrinc1}) holds $\P[\cdot\vert \xi\cap B_0\neq\emptyset]$-a.s.. The conclusion then follows by using {\it e.g. }\cite[Lemma 7.14]{CFP} wich state that a $\P[\cdot\vert \xi\cap B_0\neq\emptyset]$-almost sure event is also $\P_0$-almost sure.
\subsection{From Theorem \ref{TheoSousLine} to Theorem \ref{QIP}}\label{ProofTHPrincQIP}
Thanks to \cite[Lemma B.2]{CFP}, Theorem \ref{QIP} is a direct consequence of the following one which is a rewritting of Theorem \ref{QIP} under $\P_0$.
\begin{theo}\label{QIPPalm}
Under the assumptions of Theorem \ref{QIP}, for $\mathcal{P}_0-$a.e. $\xi^0$, under $P^{\xi^0}_0$, the rescaled process $(X^\varepsilon_t)_{t\geq 0}=(\varepsilon X_{\varepsilon^{-2}t})_{t\geq 0}$ converges in law as $\varepsilon$ tends to $0$ to a Brownian motion with covariance matrix $\sigma^2I$ where $\sigma^2$ is positive and does not depend on $\xi$. 
\end{theo}

As mentioned in the introduction, the arguments to deduce this result from Theorem \ref{TheoSousLine} are now quite standard and we only sketch the main lines of the proof. The reader is refered to \cite[\S 3.3]{CFP}, \cite[p. 1340-1341]{BP} or \cite[\S 6.1 and \S 6.2]{BergerBiskup} for more details.

Recall that, for $\P_0$-a.e. $\xi^0$, $\varphi(\xi^0,\cdot )$ is harmonic. Hence,  $(M_n)_{n\in\NN}:=(\varphi(\xi^0, X_n))_{n\in\NN}$ is a martingale under $P^{\xi^0}_0$. By the same arguments as in \cite[pp. 108-109]{BergerBiskup}, one can show that $(M_n\cdot e_i)_{n\in\NN}$ satisfies the assumptions of the Lindeberg-Feller theorem for martingales (see \cite[Theorem (7.3), p. 414]{DurrettProbab}). It follows with the Cram\'er-Wold device (see \cite[Theorem (9.5), p. 170]{DurrettProbab}) that $t\rightarrow \varepsilon M_{\lfloor \varepsilon^{-2}t\rfloor}$ converges weakly to a Brownian motion with explicit covariance matrix proportional to the identity due to the isotropy of the point process. The diffusion coefficient does not depend on the particular realization $\xi^0$ of the point process and is positive. If it were zero, it would hold that $x=\chi(\xi^0,x)$ for $\P_0$-a.e. $\xi^0$, for all $x\in\xi^0$, which contradicts the sublinearity of the corrector. The sublinearity of the corrector also implies that $\max_{k\leq n}\Vert X_k-M_k\Vert=\max_{k\leq n}\Vert\chi(\xi^0,X_k)\Vert =o(\sqrt{n})$ in $P^{\xi^0}_0$-probability. The `discrete time version' of Theorem \ref{QIPPalm} then follows. One concludes in the continuous-time case by arguing as in \cite[p. 666]{CFP} and by showing that: 
\[\lim_{t\to\infty}\frac{N(t)}{t}=\E_0\left[\deg_{\DT(\xi^0)}(0)\right],\]
where $N(t)$ denotes the number of jumps of $(X_s)_{s\geq 0}$ up to time $t$. This also proves the relation between $\sigma^2_{\text{VSRW}}$ and $\sigma^2_{\text{DTRW}}$. One finally deduces Theorem \ref{QIP} using \cite[Lemma B.2]{CFP}.  

\section{Bounds for moments of $\deg_{\DT (\xiz)}(0)$ and $\max_{x\sim 0}\Vert x\Vert$}
The method developed in \cite[\S 2]{Zuyev} can be used to derive exponential moments for $\deg_{\DT (\xiz)}(0)$ and $\max_{x\sim 0\tiny\text{ in }\DT (\xiz)}\Vert x\Vert$ when the point process has a finite range of dependence. More precisely, we show the following lemma.

\begin{lemm}\label{MomentsExpo}
Assume that $\P_0$ is isotropic and satisfies {\bf (V')}, then there exists $\rho_1>0$ such that:
\begin{equation}\label{MomentsExpoMax}
\E_0\Big[e^{\rho\max_{x\underset{\tiny\DT}{\sim} 0}\Vert x\Vert}\Big]<\infty,\quad \forall \rho<\rho_1,
\end{equation}

Assume moreover that $\P_0$ has a finite range of dependence $m$ and satisfies {\bf (EM)}, then there exists $\rho_2>0$ such that:
\begin{equation}\label{MomentsExpoDeg}
\E_0\Big[e^{\rho\deg_{\DT (\xiz)}(0)}\Big]<\infty,\quad \forall \rho<\rho_2.
\end{equation}
\end{lemm}
\begin{dem} We use the method of \cite[\S 2]{Zuyev}. Let us recall some definitions, notations, and facts from this article.
The fundamental region of a point $x\in\xiz$ is the union of the balls centered at the vertices of its Voronoi polygon and having the nucleus $x$ on their boundaries. Let $\Gamma^0$ be the union of $2d$ open balls of radii 1 centered in points $\pm e_i$. Let $\Phi^0_1,\dots,\Phi^0_{2^d}$ be the intersection of exactly $d$ such balls.

\begin{figure}[!h]
\begin{center}
\includegraphics[scale=0.17]{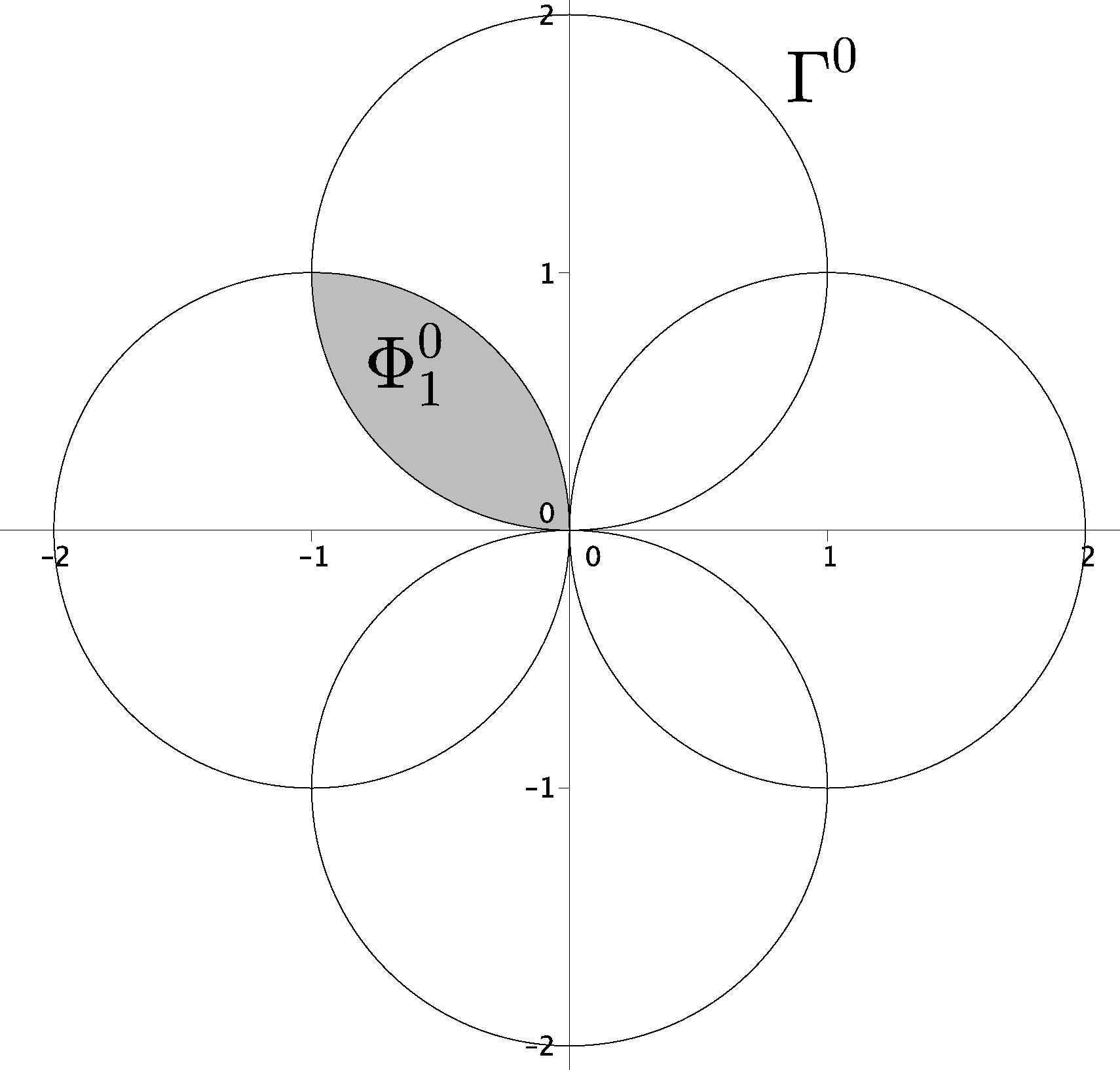}
\end{center}
\caption{$\Gamma^0$ and the lens $\Phi^0_i$.}
\end{figure}

Let us fix $\beta >1$ and consider the sequence of sets $\Gamma^n,\Phi^n_1,\dots,\Phi^n_{2^d}$ obtained by the homothetic transformation of center 0 and coefficient $\beta^n$ from $\Gamma^0,\Phi^0_1,\dots,\Phi^0_{2^d}$. The important point is that simple geometric arguments show that:

\begin{fact}\label{PropLens}
If each $d$-faced lens $\Phi^n_i, \,i=1,\dots,2^d$, contains a point of $\xiz$, then the fundamental region of the particle at 0 is fully included in $\Gamma^n$. In particular, any neighbor of 0 in $\DT (\xiz)$ is in $\Gamma^n$.
\end{fact}

Let $\A_0$ be the event $\bigcap_{i=1}^{2^d}\{\#(\xiz\cap\Phi_i^0)\neq 0\} $ and
\[\A_n:=\bigcap\limits_{i=1}^{2^d}\{\#(\xiz\cap\Phi_i^n)\neq 0\}\setminus\bigcap\limits_{i=1}^{2^d}\{\#(\xiz\cap\Phi_i^{n-1})\neq 0\}.\]

Note that the $\A_n$ are disjoint and $\P_0\big[\bigcup_{n=0}^\infty\A_n\big]=1$. Thanks to Fact \ref{PropLens}, one has:
\begin{align}\label{EqMomExpo1}
\E_0\Big[e^{\rho\max_{x\underset{\tiny\DT}{\sim} 0}\Vert x\Vert}\Big]&=\sum_{n=0}^\infty\E_0\Big[e^{\rho\max_{x\underset{\tiny\DT}{\sim} 0}\Vert x\Vert}\mathbf{1}_{\A_n}\Big]\nonumber\\
&\leq \sum_{n=0}^\infty e^{2\rho\beta^n}\P_0\big[\A_n\big]\nonumber\\
&\leq e^{2\rho}+\sum_{n=1}^\infty e^{2\rho\beta^n}\P_0\Big[\bigcup\limits_{i=1}^{2^d}\{\#(\xiz\cap\Phi_i^{n-1})= 0\}\Big]\nonumber\\
&\leq e^{2\rho}+2^d\sum_{n=1}^\infty e^{2\rho\beta^n}\P_0\Big[\#(\xiz\cap\Phi_1^{n-1})= 0\Big],
\end{align} 
where we used the isotropy of the point process in the last inequality.

Now, there exists a constant $c_{29}>0$ such that $\Phi_1^{n-1}$ contains a cube $C_{n-1}$ of side $c_{29}\beta^{n-1}$. Hence, with {\bf (V')} and (\ref{EqMomExpo1}), we obtain:
\begin{align*}
\E_0\Big[e^{\rho\max_{x\underset{\tiny\DT}{\sim} 0}\Vert x\Vert}\Big]&\leq c_{30}+2^d\sum_{n=n_0}^\infty e^{2\rho\beta^n}\P_0\Big[\#(\xiz\cap C_{n-1})= 0\Big]\\
&\leq c_{30}+2^d\sum_{n=n_0}^\infty e^{(2\rho\beta-c_{31})\beta^{(n-1)d}}.
\end{align*}
The last series converges if $\rho$ is small enough and (\ref{MomentsExpoMax}) is proved. 

In the same way, thanks to Fact \ref{PropLens}, one can write:
\begin{align}\label{EqMomExpo2}
&\E_0\Big[e^{\rho\deg_{\DT(\xiz)}(0)}\Big]\nonumber\\
&\qquad=\sum_{n=0}^\infty\E_0\Big[e^{\rho\deg_{\DT(\xiz)}(0)}\mathbf{1}_{\A_n}\Big]\nonumber\\
&\qquad\leq \E_0\Big[e^{\rho\#(\xiz\cap\Gamma^0)}\Big]+\sum_{n=1}^\infty\E_0\Big[e^{\rho\#(\xiz\cap\Gamma^n)}\mathbf{1}_{\A_n}\Big]\nonumber\\
&\qquad\leq \E_0\Big[e^{\rho\#(\xiz\cap\Gamma^0)}\Big]+\sum_{n=1}^\infty\sum_{k=0}^\infty e^{\rho k}\P_0\Big[\{\#(\xiz\cap\Gamma^n)=k\}\cap\A_n\Big]\nonumber\\
&\qquad\leq \E_0\Big[e^{\rho\#(\xiz\cap\Gamma^0)}\Big]\nonumber\\
&\qquad\qquad+2^d\sum_{n=1}^\infty\sum_{k=0}^\infty e^{\rho k}\P_0\Big[\{\#(\xiz\cap\Gamma^n)=k\}\cap\{\#(\xiz\cap\Phi_1^{n-1})=0\}\Big]\nonumber\\
&\qquad\leq\E_0\Big[e^{\rho\#(\xiz\cap\Gamma^0)}\Big]\nonumber\\
&\qquad\qquad+2^d\sum_{n=1}^\infty\sum_{k=0}^\infty e^{\rho k}\P_0\Big[\{\#(\xiz\cap\Gamma^n\setminus\Phi_1^{n-1})=k\}\cap\{\#(\xiz\cap\Phi_1^{n-1})=0\}\Big].\nonumber
\end{align} 

Now, there exists a constant $c_{32}>0$ such that, if $n$ is large enough $\Phi_1^{n-1}$ contains a cube $C_{n-1}$ of side $c_{32}\beta^{n-1}$ satisfying $\d(C_{n-1},\Gamma^n\setminus\Phi_1^{n-1})>m$. Thanks to the $m$-dependence assumption on the point process and to {\bf (EM)}, it follows that:
\begin{align*}
\E_0&\Big[e^{\rho\deg_{\DT(\xiz)}(0)}\Big]\\
&\leq c_{33}+2^d\sum_{n=n_0}^\infty\sum_{k=0}^\infty e^{\rho k}\P_0\Big[\#(\xiz\cap\Gamma^n\setminus\Phi_1^{n-1})=k\Big]\P_0\Big[\#(\xiz\cap C_{n-1})=0\Big]\\
&\leq c_{33}+ 2^d\sum_{n=n_0}^\infty\E_0\Big[e^{\rho\#(\xiz\cap [-2\beta^n,2\beta^n]^d)}\Big]\P_0\Big[\#(\xiz\cap C_{n-1})=0\Big].
\end{align*}
Finally, with {\bf (V')} and {\bf (EM)}, we obtain:
\[\E_0\Big[e^{\rho\deg_{\DT(\xiz)}(0)}\Big]\leq c_{34}+ 2^dc_5\sum_{n=n_1}^\infty e^{(2^df(\rho)\beta^d-c_{35})\beta^{(n-1)d}}.\]
This concludes the proof since $f(\rho)$ goes to 0 with $\rho$ by assumption.
\end{dem}

\subsection*{Acknowledgements}
The author thanks Jean-Baptiste Bardet and Pierre Calka, his PhD advisors, for introducing him to this subject and for helpful discussions, comments and suggestions. This work was partially supported by the French ANR grant PRESAGE (ANR-11-BS02-003) and the French research group GeoSto (CNRS-GDR3477).
\bibliography{biblio}
\bibliographystyle{alpha}
\end{document}